\documentclass[runningheads]{llncs}

\usepackage[utf8]{inputenc}
\usepackage[T1]{fontenc}
\usepackage{lmodern}
\usepackage{babel}
\usepackage{lipsum}
\usepackage{graphicx}
\usepackage{amsmath}
\usepackage{amsfonts}
\usepackage{subcaption}
\usepackage{tikz}
\usetikzlibrary{positioning, calc}
\usepackage[ruled,linesnumbered]{algorithm2e}
\usepackage{accents}
\usepackage{xcolor}
\usepackage{parskip}
\usepackage{amssymb}
\usepackage[hyperfootnotes=false]{hyperref}
\usepackage{thmtools}
\usepackage{multirow}
\usepackage{alphabeta} 
\usepackage{thm-restate}

\newif\ifarxiv
\arxivtrue

\newif\ifopus
\opusfalse

\usepackage[scaled=0.95]{FiraMono}
\ifarxiv
\usepackage{charter}
\usepackage{fullpage}
\fi

\ifopus
\usepackage{zibtitlepage}
\fi
\usepackage{todonotes}

\DeclareMathOperator*{\argmin}{arg\,min}

\newcommand{\coefficients}{\ensuremath{\boldsymbol{\alpha}}\xspace}
\newcommand{\objective}{\ensuremath{\mathbf{c}}\xspace}

\newcommand{\xv}{\ensuremath{\mathbf{x}}}
\newcommand{\lpoptimal}{\ensuremath{\mathbf{x}^{LP}}\xspace}

\newcommand{\lpoptimalset}{\ensuremath{\mathcal{X}^{LP}}\xspace}
\newcommand{\incumbent}{\ensuremath{\hat{\mathbf{x}}}\xspace}
\newcommand{\lpcenface}{\ensuremath{\mathbf{x}^{F}}\xspace}
\newcommand{\lpcen}{\ensuremath{\mathbf{x}^{C}}\xspace}

\newcommand{\naturals}{\ensuremath{\mathbb{N}}\xspace}

\newcommand{\eff}[3]{\ensuremath{\mathtt{eff}(#1,#2,#3)}\xspace}
\newcommand{\dcd}[4]{\ensuremath{\mathtt{dcd}(#1,#2,#3,#4)}\xspace}
\newcommand{\avgeff}[3]{\ensuremath{\mathtt{avgeff}(#1,#2,#3)}\xspace}
\newcommand{\mineff}[3]{\ensuremath{\mathtt{mineff}(#1,#2,#3)}\xspace}
\newcommand{\expimprov}[4]{\ensuremath{\mathtt{exp}\text{-}\mathtt{improv}(#1,#2,#3,#4)}\xspace}
\newcommand{\density}[1]{\ensuremath{\mathtt{dens}(#1)}\xspace}
\newcommand{\relativedensity}[1]{\ensuremath{\mathtt{reldens}(#1)}\xspace}

\newcommand{\noargdcd}{\ensuremath{\mathtt{dcd}}\xspace}
\newcommand{\noargeff}{\ensuremath{\mathtt{eff}}\xspace}
\newcommand{\noargavgeff}{\ensuremath{\mathtt{avgeff}}\xspace}
\newcommand{\noargmineff}{\ensuremath{\mathtt{mineff}}\xspace}
\newcommand{\noarganaldcd}{\ensuremath{\mathtt{a}\text{-}\mathtt{dcd}}\xspace}
\newcommand{\noarganaleff}{\ensuremath{\mathtt{a}\text{-}\mathtt{eff}}\xspace}
\newcommand{\noargapproxanaldcd}{\ensuremath{\mathtt{app}\text{-}\mathtt{a}\text{-}\mathtt{dcd}}\xspace}
\newcommand{\noargexpimprov}{\ensuremath{\mathtt{exp}\text{-}\mathtt{improv}}\xspace}
\newcommand{\noargdensity}{\ensuremath{\mathtt{dens}}\xspace}
\newcommand{\noargrelativedensity}{\ensuremath{\mathtt{reldens}}\xspace}

\hypersetup{
    colorlinks=true,
    linkcolor=blue,
    filecolor=magenta,      
    urlcolor=cyan,
    citecolor=blue
    }

\def\hat{\mathaccent "705E\relax}

\title{Cutting Plane Selection with Analytic Centers and Multiregression}

\date{August 2022}

\begin{document}

\ifarxiv
\author{Mark Turner\inst{1,2}\and
Timo Berthold\inst{1,3}\and\\
Mathieu Besançon\inst{2}\and
Thorsten~Koch\inst{1,2}
}
\else
\author{Mark Turner\inst{1,2}\orcidID{0000-0001-7270-1496} \and
Timo Berthold\inst{1,3}\orcidID{0000-0002-6320-8154}\and\\
Mathieu Besançon\inst{2}\orcidID{0000-0002-6284-3033}\and
Thorsten~Koch\inst{1,2}\orcidID{0000-0002-1967-0077}
}
\fi

\authorrunning{M. Turner et al.}

\institute{Institute of Mathematics, Technische Universit{\"a}t Berlin, Germany \and
Zuse Institute Berlin, Germany\and
Fair Isaac Deutschland GmbH, Berlin, Germany\\
\email{\{turner, besancon, koch\}@zib.de}, \email{timoberthold@fico.com}
}

\ifopus
\ZTPAuthor{\ZTPHasOrcid{Mark Turner}{0000-0001-7270-1496},
\ZTPHasOrcid{Timo Berthold }{0000-0002-6320-8154},
\ZTPHasOrcid{Mathieu Besançon}{0000-0002-6284-3033}, \\
\ZTPHasOrcid{Thorsten Koch }{0000-0002-1967-0077}}
\ZTPTitle{Cutting Plane Selection with Analytic Centers and Multiregression}
\ZTPNumber{22-28}
\ZTPMonth{December}
\ZTPYear{2022}
\fi

\ifopus
\fi

\maketitle

\begin{abstract}
Cutting planes are a crucial component of state-of-the-art mixed-integer programming solvers, with the choice of which subset of cuts to add being vital for solver performance.
We propose new distance-based measures to qualify the value of a cut by quantifying the extent to which it separates relevant parts of the relaxed feasible set.
For this purpose, we use the analytic centers of the relaxation polytope or of its optimal face, as well as alternative optimal solutions of the linear programming relaxation. 
We assess the impact of the choice of distance measure on root node performance and throughout the whole branch-and-bound tree, comparing our measures against those prevalent in the literature.
Finally, by a multi-output regression, we predict the relative performance of each measure, using static features readily available before the separation process.
Our results indicate that analytic center-based methods help to significantly reduce the number of branch-and-bound nodes needed to explore the search space and that our multiregression approach can further improve on any individual method.
\end{abstract} 

\section{Introduction}
\label{sec:introduction}

Branch-and-cut is the algorithm at the core of most Mixed-Integer Programming (MIP) solvers.
A key component of branch-and-cut is the resolution of Linear Programming (LP) relaxations of the original problem over partitions of the variable domains.
\emph{Cutting planes} -- or cuts --
tighten those relaxations around integer-feasible points. Given a MIP:
\begin{align}
    \underset{\mathbf{x}}{\text{argmin}}\{\mathbf{c}^{\intercal}\mathbf{x} \;\; | \;\; \mathbf{A}\mathbf{x} \leq \mathbf{b}, \;\; \mathbf{l} \leq \mathbf{x} \leq \mathbf{u}, \;\; \mathbf{x} \in \mathbb{Z}^{|\mathcal{J}|} \times \mathbb{R}^{n - |\mathcal{J}|} \} \label{eq:mip}
    \tag{P}
\end{align}
a cut is an inequality $\coefficients^{\intercal} \xv \leq \beta$ that is violated by at least one solution of the LP relaxation but that does not increase the optimal value of the problem when added i.e.~it is valid for \eqref{eq:mip}.
Thereby, the inequality added as a constraint to \eqref{eq:mip} tightens the relaxation, potentially increasing the relaxation's optimal value.
The use of cutting planes is one of the crucial aspects to solving MIPs efficiently~\cite{achterberg2013mixed}.
A well-designed cutting plane separation procedure often helps to reduce the branch-and-bound tree size while accelerating the overall solving process.

In MIP solvers, two key algorithms related to cuts are their generation and their selection.
Cut generation is the problem of computing a set of cuts that tighten the relaxation and separate the current continuous relaxation solution from the feasible MIP solutions.
Modern MIP solvers implement various general-purpose and specialised cutting plane generation algorithms.
Since the generation of cuts is, in general, far less expensive than solving the LP relaxation, many cuts are generated from the same relaxation. 
The cut selection algorithm takes the set of all cut candidates generated so far
and selects a subset that is actually added to the LP relaxation.
This two-step process of generation and selection constitutes a single \emph{separation round}.

At the root node, the MIP solver interleaves separation rounds with solving the enhanced LP relaxation until the branch-and-bound search is started. At other search tree nodes, the solver often only performs a limited cut loop, if at all.
We focus in this paper on globally-valid cuts, i.e.~cuts that are valid for the original
problem, as opposed to \emph{local cuts}, which are generated with additional local bounds at a node.

Cut selection is a classical trade-off problem: too little cutting leads to large enumeration trees; too much cutting to a small node throughput and numerical instability.
However, carefully selected cuts can help improve both the dual and the primal bound simultaneously by bringing the relaxation closer to the convex hull of feasible solutions. Since proximity to the convex hull is inherently hard to measure, cut selection methods often try to approximate it by various cheap measures, e.g., efficacy.



\emph{Efficacy}, or \emph{cutoff distance}, is used in commercial MIP solvers as one of the main criteria for whether to add a cut or not~\cite{achterberg2007constraint,scip8}\footnote{Confirmed as a main criterion in FICO Xpress 8.14}. Efficacy measures the shortest distance between the LP optimal solution \lpoptimal and the cut hyperplane $\coefficients^{\intercal} \mathbf{x} \leq \beta$. The function \noargeff that maps a cut, and the LP solution \lpoptimal to the efficacy is defined as:
\begin{align*}
    \eff{\coefficients}{\beta}{\lpoptimal} := \frac{\coefficients^{\intercal} \lpoptimal  - \beta}{\|\coefficients\| }
\end{align*}

Introduced in \cite{scip6}, \emph{directed cutoff distance} is the signed distance between the LP solution \lpoptimal and the cut hyperplane in the direction of a primal solution, \incumbent. The measure has the property that the directed projection of \lpoptimal onto the cut hyperplane is inside of the feasible region, and by using the best primal solution available, aims to cut in the direction of the optimal solution. We define the directed cutoff distance, with the function \noargdcd as follows:
\begin{align*}
    \dcd{\coefficients}{\beta}{\lpoptimal}{\incumbent} := \frac{\coefficients^{\intercal} \lpoptimal  - \beta}{|\coefficients^{\intercal} \mathbf{y}|}, \; \text{where} \quad \mathbf{y} = \frac{\incumbent - \lpoptimal}{||\incumbent - \lpoptimal||}
\end{align*}

The last measure we consider, although not based on a distance, is \emph{expected improvement}, see \cite{wesselmann2012implementing},
which corresponds to the difference in objective between \lpoptimal and its orthogonal projection onto the hyperplane of a cut. We denote the measure as \noargexpimprov, and define it as:
\begin{align*}
    \expimprov{\coefficients}{\beta}{\objective}{\lpoptimal} :=& \; ||\objective|| \; \cdot \; \frac{\coefficients^{\intercal}\objective}{||\coefficients||||\objective||} \; \cdot \; \eff{\coefficients}{\beta}{\lpoptimal}
\end{align*}

In this paper, we propose new distance-based measures for the quality of cuts.
These measures are designed to retain soundness properties in common cases that hinder the applicability, and to reduce the size of the search space, without focusing on runtime improvement.
To this end, we perform extensive computational experiments
to analyse the effects of our newly introduced measures and those that exist in the literature.
Unlike previous beliefs, the choice of distance measure significantly impacts solver performance both in time and number of nodes, with different measures performing better on different groups of instances.
Motivated by this observation, we design a multi-output regression model, which predicts the relative performance of each measure using static features readily available before the separation process. 
\ifarxiv
Finally, in Appendix \ref{sec:density}, we complement our investigation of cut selection with an analysis of the impact of dense cutting planes on solver performance.
\fi
The scope of this paper is to establish a model that aims at reducing search space size rather than runtime.

\section{Related Work}

There is a prevailing sentiment in the MIP community, supported by a set of older computational studies, see \cite{achterberg2007constraint,zerohalf,wesselmann2012implementing}, that inexpensive heuristics are sufficient for cut selection.
Specifically, these studies suggest that cheap ranking metrics, predominantly efficacy, are sufficiently effective when combined with filtering mechanisms that ensure no two overly parallel cuts are added. The studies \cite{achterberg2007constraint,wesselmann2012implementing} argue that a weighted sum of different metrics is most effective for ranking cuts as opposed to any single metric.

Recently, the research focus for cuts in MIP has been on using deep learning to either calculate scores directly from a set of measures or to predict parameter values in scoring functions, see \cite{baltean2019scoring,tang2020reinforcement,huang2021learning,balcan2021sample,turner2022adaptive,paulus2022learning}. In \cite{baltean2019scoring} and \cite{paulus2022learning} a neural network is trained to predict the objective value improvement of cuts when added. In \cite{tang2020reinforcement}, a neural network is trained using evolutionary strategies to select Gomory cuts. A neural network is also trained in \cite{huang2021learning}, this time using multiple instance learning, to map a set of aggregate cut features to a scoring function. In \cite{balcan2021sample}, the cut selection parameter space is shown to be partitionable into a set of regions, such that all parameter choices within a region select the same set of cuts.
Finally, for cut scores based on weighted criteria, \cite{turner2022adaptive} provides an illustrative example of worst-case scenario for parameter grid search in the cut selection parameter space, and phrases learning cut selection parameters as a reinforcement learning problem.

Closest to our work are papers that introduced cut selection measures other than efficacy, which were still based on a notion of measuring distances.
\emph{Directed cutoff distance} was introduced in SCIP 6.0~\cite{scip6}, other measures such as \emph{rotated cutoff distance} and \emph{distance with bounds} were explored in \cite{wesselmann2012implementing}, and \emph{depth} was introduced in \cite{depth}. Measures based on non-distance arguments are also prevalent in cut selection, albeit often as smaller weighted complements to a distance measure, see \emph{objective parallelism} and \emph{integer support} in \cite{achterberg2007constraint,wesselmann2012implementing}, and enumeration of lattice points in \cite{latticepoints}.

In the presented distance measures of Section \ref{sec:introduction}, an LP-optimal solution is used as a reference point.
Note, however, that this optimum is not necessarily unique, with the set of minimisers frequently being a higher-dimensional face of the LP-feasible region due to dual degeneracy~\cite{mipdegeneracy}.
The optimal solution used in these degenerate cases inevitably has a large impact on cut generation~\cite{lexicographicgomory}, and therefore cut selection. This has been noted in previous research, with work such as \cite{multiplelpcutting} using multiple LP solutions from different LP random seeds to generate different sets of cuts. Additionally, the patent \cite{cplexpatent}, proposes using a second LP optimal solution at the same cutting round to filter cuts derived from the original LP solution. They provide an example, where a second LP solution that prioritises integrality is found, which then filters all cuts that do not separate it. Dual degeneracy is one major aspect motivating some of the newly-proposed measures in our work.

Finally, analytic centers have been used in other aspects of MIP solving, namely, presolving~\cite{berthold2018four}, cut generation~\cite{fischetti2009yoyo}, branching~\cite{berthold2018four}, and heuristics~\cite{BAENA2011310,naoum2014recursive}, motivating the measures introduced in this paper.

\section{Contributions and methodology}

The contribution of the present paper is threefold.
First, we introduce new distance-derived cut quality measures, the most important of which utilise analytic centers, and analyse properties of interest for these measures in cases of dual degeneracy and infeasible projections. Secondly, we present an extensive set of computational experiments on the effectiveness of our new measures and those commonly found in practice,
showing that the choice of cut selection measure does have a strong influence on root-node and tree-wide performance.
Thirdly, we introduce a multi-output regression model that predicts a ranking of distance measures per instance from a set of root node features.

\subsection{Analytic Center-Based Methods}

We propose two new methods for measuring cut quality: \emph{analytic efficacy}, and \emph{analytic directed cutoff distance}. They are based on the analytic center of the polytope and of the optimal face, respectively. For a given bounded MIP formulation \eqref{eq:mip}, the analytic center is unique and in the relative interior of the feasible region.

When a constraint $\mathbf{a}^{\intercal} \mathbf{x} \leq b$ is tight for any feasible solution
or in the presence of equality constraints,
the analytic center is not well-defined due to a log-barrier term being $+\infty$.
In practice, algorithms relax all log-barriers with a fixed slack constant
as long as constraints are imposed on solutions.

While the analytic center is invariant under affine transformations of Problem~\eqref{eq:mip}, it can change with reformulations. E.g., the analytic center can be shifted by the presence
of redundant inequalities. For our work, we assume that the formulation has already been presolved by the MIP solver.
Presolving includes tightening both variable bounds and constraints, and removing both redundant constraints and variables. It thereby limits the extent of such edge cases.


\subsubsection{Analytic efficacy}\label{subsec:ae}

As opposed to using the \lpoptimal extreme point returned by the MIP solver for \noargeff calculations, we propose a new measure that uses the analytic center of the optimal face of the LP, \lpcenface. We define \lpcenface by the following using the notation from Problem~\eqref{eq:mip}, where $\mathbf{A}_{i}$, $\mathbf{b}_{i}$, $\mathbf{u}_{i}$, and $\mathbf{l}_{i}$ are the $i$-th row or $i$-th entry of their respective matrix or vector:
\begin{align*}
    \lpcenface &:= \underset{\mathbf{x}\; : \; \mathbf{c}^{\intercal}\mathbf{x} = \mathbf{c}^{\intercal}\lpoptimal}{\text{argmin}}\{-\sum_{i=1}^{m} \log ( \mathbf{b}_{i} - \mathbf{A}_{i} \mathbf{x} ) - (\sum_{i=1}^{n} \log(\mathbf{x}_{i} - \mathbf{l}_{i}) + \log(\mathbf{u}_{i} - \mathbf{x}_{i})) \}
\end{align*}
In practice, we find \lpcenface using the barrier (or interior point) algorithm on the LP relaxation without crossover, see \cite{berthold2018four}. This algorithm is available in all modern MIP solvers, and is often run concurrently to the simplex at the root node, making our cut selection algorithms of practical interest for MIP solving.

The obtained center \lpcenface is different from \lpoptimal only if the current problem
presents dual degeneracy.
The purpose of evaluating cuts with respect to how much they separate the center of the optimal face is to favour those which cut a greater part or potentially all of the optimal face, thereby more likely favouring an improvement in the dual bound. Compare Figures~\ref{fig:efficacy_visualisation} and~\ref{fig:analytic_efficacy_visualisation} for the intuition behind analytic efficacy.

\subsubsection{Analytic directed cutoff distance}

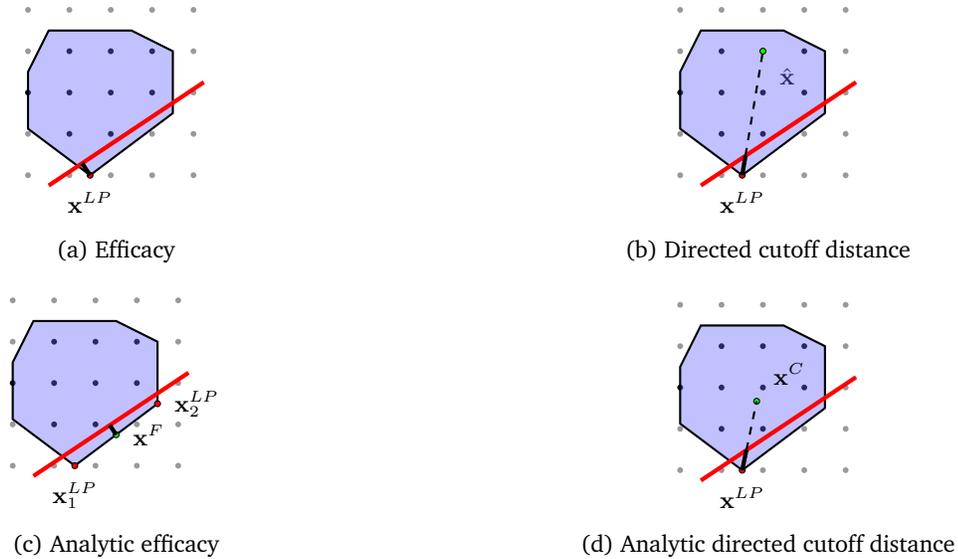
\begin{figure*}[h]
\centering
\begin{subfigure}[b]{0.475\textwidth}
    \centering
    \begin{tikzpicture}[scale=0.55]
    \foreach \x/\y in {1/0, 1/1, 1/3, 1/4, 2/0, 2/4, 3/0, 3/4, 4/0, 4/1, 4/4, 5/0, 5/1, 5/2, 5/3, 5/4}
    {
    \fill[opacity=0.4] (\x,\y) circle (2pt);
    }
    \foreach \x/\y in {2/3, 3/3, 4/3, 1/2, 2/2, 3/2, 4/2, 2/1, 3/1}
    {
    \fill[opacity=1] (\x,\y) circle (2pt);
    }
    \node (a) at (1,1.125) {};
    \node (b) at (1,2.5) {};
    \node (c) at (1.5,3.5) {};
    \node (d) at (3.5,3.5) {};
    \node (e) at (4.5,3) {};
    \node (f) at (4.5,1.5) {};
    \node[label=below:{$\lpoptimal$}] (g) at (2.5,0) {};
    
    \filldraw[thick,fill=blue!60,fill opacity=0.4] (a.center) -- (b.center) -- (c.center) -- (d.center)  -- (e.center) -- (f.center) -- (g.center) -- cycle;
    
    \node (lcut) at (1.5,-0.25) {};
    \node (rcut) at (5.25,2.25) {};
    \draw[ultra thick, red] (lcut.center) -- (rcut.center);
    
    \fill[draw=black,fill=red] (2.5,0) circle (2pt);
    
    \draw[ultra thick] (30/13,15/52) -- (2.5, 0);
    
    
    
    
    
    \end{tikzpicture}
\caption{Efficacy}
\label{fig:efficacy_visualisation}
\end{subfigure}
\hfill
\begin{subfigure}[b]{0.475\textwidth}
\centering
    \begin{tikzpicture}[scale=0.55]
    \foreach \x/\y in {1/0, 1/1, 1/3, 1/4, 2/0, 2/4, 3/0, 3/4, 4/0, 4/1, 4/4, 5/0, 5/1, 5/2, 5/3, 5/4}
    {
    \fill[opacity=0.4] (\x,\y) circle (2pt);
    }
    \foreach \x/\y in {2/3, 3/3, 4/3, 1/2, 2/2, 3/2, 4/2, 2/1, 3/1}
    {
    \fill[opacity=1] (\x,\y) circle (2pt);
    }
    
    \node[label=below right:{$\incumbent$}] (primal) at (3,3) {};
    \node (a) at (1,1.125) {};
    \node (b) at (1,2.5) {};
    \node (c) at (1.5,3.5) {};
    \node (d) at (3.5,3.5) {};
    \node (e) at (4.5,3) {};
    \node (f) at (4.5,1.5) {};
    \node[label=below:{$\lpoptimal$}] (g) at (2.5,0) {};
    
    \filldraw[thick,fill=blue!60,fill opacity=0.4] (a.center) -- (b.center) -- (c.center) -- (d.center)  -- (e.center) -- (f.center) -- (g.center) -- cycle;
    
    \node (lcut) at (1.5,-0.25) {};
    \node (rcut) at (5.25,2.25) {};
    \draw[ultra thick, red] (lcut.center) -- (rcut.center);
    
    \fill[draw=black,fill=red] (2.5,0) circle (2pt);
    \fill[draw=black,fill=green] (3,3) circle (2pt);
    
    \draw[ultra thick] (2.5,0) -- (165/64,15/32);
    \draw[thick,dashed] (165/64,15/32) -- (3,3);
    \end{tikzpicture}
\caption{Directed cutoff distance}
\label{fig:directed_cutoff_visualisation}
\end{subfigure}
\vskip\baselineskip
\begin{subfigure}[b]{0.475\textwidth}
    \centering
    \begin{tikzpicture}[scale=0.55]
   \foreach \x/\y in {1/0, 1/1, 1/3, 1/4, 2/0, 2/4, 3/0, 3/4, 4/0, 4/1, 4/4, 5/0, 5/1, 5/2, 5/3, 5/4}
    {
    \fill[opacity=0.4] (\x,\y) circle (2pt);
    }
    \foreach \x/\y in {2/3, 3/3, 4/3, 1/2, 2/2, 3/2, 4/2, 2/1, 3/1}
    {
    \fill[opacity=1] (\x,\y) circle (2pt);
    }
    \node (a) at (1,1.125) {};
    \node (b) at (1,2.5) {};
    \node (c) at (1.5,3.5) {};
    \node (d) at (3.5,3.5) {};
    \node (e) at (4.5,3) {};
    \node (f) at (4.5,1.5) {};
    \node[label=below:{$\lpoptimal_1$}] (g) at (2.5,0) {};
    
    \filldraw[thick,fill=blue!60,fill opacity=0.4] (a.center) -- (b.center) -- (c.center) -- (d.center)  -- (e.center) -- (f.center) -- (g.center) -- cycle;
    
    \node (lcut) at (1.5,-0.25) {};
    \node (rcut) at (5.25,2.25) {};
    \draw[ultra thick, red] (lcut.center) -- (rcut.center);
    
    \fill[draw=black,fill=red] (2.5,0) circle (2pt);
    \fill[draw=black,fill=green] (3.5,0.75) circle (2pt);
    \fill[draw=black,fill=red] (4.5,1.5) circle (2pt);
    
    \node[label=right:{\lpcenface}] (facecenter) at (3.5,0.75) {};
    \node[label=right:{$\lpoptimal_2$}] (facecenter) at (4.5,1.5) {};
    
    \draw[ultra thick] (3.5,0.75) -- (87/26, 51/52);
    
    \end{tikzpicture}
\caption{Analytic efficacy}
\label{fig:analytic_efficacy_visualisation}
\end{subfigure}
\hfill
\begin{subfigure}[b]{0.475\textwidth}
    \centering
    \begin{tikzpicture}[scale=0.55]
    \foreach \x/\y in {1/0, 1/1, 1/3, 1/4, 2/0, 2/4, 3/0, 3/4, 4/0, 4/1, 4/4, 5/0, 5/1, 5/2, 5/3, 5/4}
    {
    \fill[opacity=0.4] (\x,\y) circle (2pt);
    }
    \foreach \x/\y in {2/3, 3/3, 4/3, 1/2, 2/2, 3/2, 4/2, 2/1, 3/1}
    {
    \fill[opacity=1] (\x,\y) circle (2pt);
    }
    \node (a) at (1,1.125) {};
    \node (b) at (1,2.5) {};
    \node (c) at (1.5,3.5) {};
    \node (d) at (3.5,3.5) {};
    \node (e) at (4.5,3) {};
    \node (f) at (4.5,1.5) {};
    \node[label=below:{$\lpoptimal$}] (g) at (2.5,0) {};
    
    \filldraw[thick,fill=blue!60,fill opacity=0.4] (a.center) -- (b.center) -- (c.center) -- (d.center)  -- (e.center) -- (f.center) -- (g.center) -- cycle;
    
    \node (lcut) at (1.5,-0.25) {};
    \node (rcut) at (5.25,2.25) {};
    \draw[ultra thick, red] (lcut.center) -- (rcut.center);
    
    \fill[draw=black,fill=red] (2.5,0) circle (2pt);
    \fill[draw=black,fill=green] (2.849,1.667) circle (2pt);
    \node[label=above right:{\lpcen}] (anycenter) at (2.849,1.667) {};
    
    \draw[ultra thick] (2.5,0) -- (2.601,0.4842);
    \draw[thick, dashed] (2.601,0.4842) -- (2.849,1.667);
    \end{tikzpicture}
\caption{Analytic directed cutoff distance}
\label{fig:analytic_directed_cutoff_visualisation}
\end{subfigure}
\caption{A visualisation of distance measures. Note that in Figure~\ref{fig:analytic_efficacy_visualisation} we see that there are two alternative optimal vertices.
}
\end{figure*}

This measure is inspired by directed cutoff distance, and uses the analytic center of the feasible region, \lpcen, as opposed to the best incumbent solution \incumbent.
Using the same notation as in Section \ref{subsec:ae}, we define the analytic center of the feasible region, \lpcen, as:
\begin{align*}
    \lpcen &:= \underset{\mathbf{x}}{\text{argmin}}\{-\sum_{i=1}^{m} \log ( \mathbf{b}_{i} - \mathbf{A}_{i} \mathbf{x} ) - (\sum_{i=1}^{n} \log(\mathbf{x}_{i} - \mathbf{l}_{i}) + \log(\mathbf{u}_{i} - \mathbf{x}_{i})) \}
\end{align*}
\lpcen is often interpreted as the ``central-most'' point of the polytope
and can be computed efficiently by dropping the objective function $\mathbf{c}^{\intercal} \xv$ and solving the resulting LP relaxation using the barrier algorithm without crossover.
The motivation is that the incumbent is not necessarily representative of the feasible set of the MIP since it can be any point of the feasible set generated by heuristics.
Furthermore, it introduces an additional source of variability in the cut selection since the best primal solution is prone to updates, particularly during root node cutting when many primal heuristics are employed. By contrast, \lpcen is unique and deterministically determined for the LP relaxation. Compare Figures~\ref{fig:directed_cutoff_visualisation} and~\ref{fig:analytic_directed_cutoff_visualisation} for the intuition behind analytic directed cutoff distance.

For computational efficiency, we further introduce \emph{approximate analytic directed cutoff distance}, which re-uses the analytic center, \lpcen, from the previous separation round provided it is still LP-feasible. This is motivated by the intuition that the analytic center, as the ``central-most'' point, is rarely separated and remains close to the new analytic center after cuts have been added.

\subsection{Multiple LP solutions}
\label{subsec:multiple_lp}

We also introduce two distance measures, which mitigate cases of dual degeneracy and do not rely on analytic centers, but rather rely on multiple LP-optimal vertices.
Let \lpoptimalset be a set of LP optimal solutions. We propose the measures \emph{average efficacy}, denoted \noargavgeff, and \emph{minimum efficacy}, denoted \noargmineff, which respectively take the average and minimum efficacy over all LP solutions in \lpoptimalset.
\begin{align*}
    \avgeff{\coefficients}{\beta}{\lpoptimalset} &:= \sum_{\lpoptimal \in \lpoptimalset} \frac{\eff{\coefficients}{\beta}{\lpoptimal}}{|\lpoptimalset|} \\
    \mineff{\coefficients}{\beta}{\lpoptimalset} &:= \text{min}\{ \eff{\coefficients}{\beta}{\lpoptimal} \;\; | \;\; \lpoptimal \in \lpoptimalset \}
\end{align*}

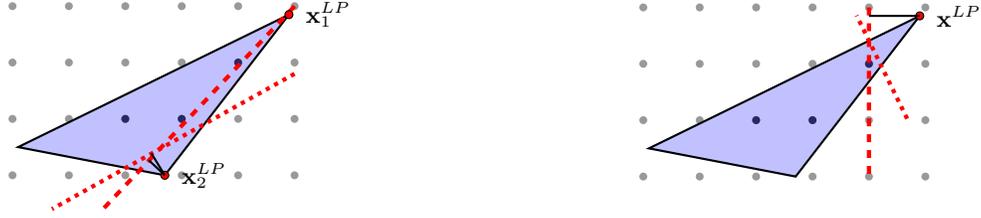
\begin{figure}[htb]
\begin{minipage}[t]{.492\textwidth}
    \centering
    \begin{tikzpicture}[rotate=-90,scale=0.75]
    \foreach \x/\y in {3/1, 3/2, 3/3, 3/4, 3/5, 3/6, 4/1, 4/2, 4/3, 4/4, 4/5, 4/6, 5/1, 5/2, 5/3, 5/4, 5/5, 5/6, 6/1, 6/2, 6/3, 6/4, 6/5, 6/6}
    {
    \fill[opacity=0.4] (\x,\y) circle (2pt);
    }
    \foreach \x/\y in {4/5, 5/3, 5/4}
    {
    \fill[opacity=1] (\x,\y) circle (2pt);
    }
    \node[label=right:{$\lpoptimal_1$}] (a) at (3.15, 5.9) {};
    \node[label=right:{$\lpoptimal_2$}] (b) at (6, 3.7) {};
    \node (c) at (5.5, 1.1) {};
    
    \filldraw[thick,fill=blue!60,fill opacity=0.4] (a.center) -- (b.center) -- (c.center) -- cycle;
    
    \draw[ultra thick, dashed, red] (3,6) -- (6.6771, 2.5313);
    \draw[ultra thick, dotted, red] (4.2,6) -- (6.6,1.7);
    
    \fill[draw=black,fill=red] (6,3.7) circle (2pt);
    \fill[draw=black,fill=red] (3.15,5.9) circle (2pt);
    
    \draw[thick] (6, 3.7) -- (5.611, 3.462);
    \draw[thick] (6, 3.7) -- (5.7355, 3.4196);
    
    
    \end{tikzpicture}
\subcaption{The dashed cut removes the entire optimal face, and thereby is likely preferable. Efficacy however scores the dotted cut higher (black lines). Analytic, minimum, and average efficacy prefer the dashed cut.}\label{fig:eff_dual_degeneracy}
\end{minipage}
\hfill
\begin{minipage}[t]{.492\textwidth}
    \centering
    \begin{tikzpicture}[rotate=-90,scale=0.75]
    \foreach \x/\y in {3/1, 3/2, 3/3, 3/4, 3/5, 3/6, 4/1, 4/2, 4/3, 4/4, 4/5, 4/6, 5/1, 5/2, 5/3, 5/4, 5/5, 5/6, 6/1, 6/2, 6/3, 6/4, 6/5, 6/6}
    {
    \fill[opacity=0.4] (\x,\y) circle (2pt);
    }
    \foreach \x/\y in {4/5, 5/3, 5/4}
    {
    \fill[opacity=1] (\x,\y) circle (2pt);
    }
    \node[label=right:{$\lpoptimal$}] (a) at (3.15, 5.9) {};
    
    \filldraw[thick,fill=blue!60,fill opacity=0.4] (a.center) -- (b.center) -- (c.center) -- cycle;
    
    \fill[draw=black,fill=red] (3.15,5.9) circle (2pt);
    
    \draw[ultra thick, dashed, red] (3,5) -- (6,5);
    \draw[ultra thick, dotted, red] (4.9742,5.677) -- (3.1426,4.7807);
    \draw[thick] (3.15,5.9) -- (3.15,5);
    
    \draw (6.6771, 2.5313) {};
    \draw (6.6,1.7) {};
    
    \end{tikzpicture}
    \subcaption{The projection from \lpoptimal onto the dashed cut is LP-infeasible, limiting its usefulness as a measure. The dotted cut is dominated but would still be selected by efficacy.}\label{fig:infeasible_projection}
\end{minipage}
\caption{Two cases showing the limitation of efficacy, in the presence of dual degeneracy \ref{fig:eff_dual_degeneracy}, and infeasible projection \ref{fig:infeasible_projection}. The blue polytope represents the LP feasible region, the black dots integer solutions, the red lines are proposed cuts, and \lpoptimal represents a LP optimal solution. In both cases, the dotted cut has a higher efficacy even though it is likely not preferable.}
\end{figure}

\subsection{Properties and limitations of the distance measures}
\label{subsec:properties}

In this section, we highlight situations in which the standard interpretations of some distance measures are limited;
we introduce dominance consistency, a soundness property for distance measures, and establish cases under which it holds and under which it does not. Examples of dominance consistency also offer valuable insight into the geometry of different measures.

Efficacy, similarly expected improvement, uses the orthogonal projection of \lpoptimal onto the cut to measure distance. Unlike measures such as directed cutoff distance, the projection point might not be LP feasible, potentially making efficacy non-representative of the strength of the cut. A larger efficacy does not necessarily correspond to a larger part of the polyhedron being cut off nor to a better improvement in dual bound, see Figure~\ref{fig:eff_dual_degeneracy}. Note that minimum and analytic efficacy would not assign the dotted cut a positive score. Efficacy would prefer the dotted cut even though it does not improve the dual bound while the dashed cut does. Figure~\ref{fig:infeasible_projection} visualises a basic example of when the orthogonal projection used for efficacy is LP-infeasible and thus not necessarily a good proxy for cut quality. Analytic efficacy, minimum efficacy, and average efficacy help overcome some limitations of efficacy, namely the dependence on the vertex returned by the LP solve.
They are, however, equivalent when the current relaxation is not dual-degenerate since they then compute the distance to the cut using the same, unique optimal vertex.

Directed cutoff distance heavily depends on the incumbent solution, which is typically obtained from a primal heuristic.\footnote{At the root node, in particular, the incumbent -- if existing -- will always come from a heuristic, otherwise there would be no more cut rounds.} This primal solution may be sub-optimal and near a corner of the feasible region, biasing cuts in an unfavourable direction. Additionally, the primal solution may be LP-infeasible for local relaxations of the branch-and-bound tree, reducing the applicability of directed cutoff distance.

Ultimately, one would like a distance measure to be a surrogate for solver efficiency, in our case: the size of the search space in terms of branch and bound nodes. Such a measure is impossible to quantify however, with solution fractionality being insufficient, and the closest analogue, `strong cutting', which measures the dual bound improvement of an added cut, see \cite{paulus2022learning}, being computationally intractable. The holy grail of cut selection is to either identify a computationally tractable measure that is a reasonably good proxy for solver performance, or to learn to adaptively select a suitable measure based on the input instance.

In the following, we will formalise these considerations. Therefore, we first recall the basic concepts of cut dominance and feasible rays.
A cut $(A) = (\coefficients_A,\beta_A)$ \emph{dominates} another cut $(B) = (\coefficients_B,\beta_B)$ if all points of the polytope cut by $(A)$ are cut by $(B)$ and there exists a point cut by $(A)$ not cut by $(B)$.
We highlight that this definition of dominance is more general than, e.g., \cite[Definition 9.2.1]{wolsey2020integer} since it only requires dominance to hold in
the polytope and not in the whole space or the positive orthant.
We define a \emph{feasible ray} as a ray $\mathbf{r}$ starting from an LP-feasible point $\xv$ for which there exists $\lambda > 0$ such that $\xv + \lambda \mathbf{r}$ is LP-feasible. Such a ray always exists by convexity if the polyhedron is not a single point.

\begin{definition}[Dominance consistency]
Given a MIP \eqref{eq:mip} and a relaxation point to cut off \xv, a distance measure noted $d(\xv, \coefficients_X,\beta_X)$
is \emph{dominance-consistent} w.r.t.~a set of cuts
iff for any cut $(A)$ and $(B)$ in the set, $d(\xv, \coefficients_A,\beta_A) > d(\xv, \coefficients_B,\beta_B)$ implies that $(A)$ is not dominated by $(B)$.
\end{definition}


Note that with this definition, for a given MIP and LP solution \xv, one set of cuts (e.g., from one separation round) can imply dominance consistency for a measure while another set of cuts (e.g., from another separation round) might not. We will see that for some measures, dominance consistency applies for all sets of cuts, all relaxation points, and all MIPs.
Dominance consistency is a desirable property for a cut selection measure
since a fully dominated cut will systematically be inferior to the dominating cut from the cut strength perspective.\footnote{Note that we disregard other cut properties such as density, numerical stability, or orthogonality here.}

\begin{proposition}[Consistency of Euclidean distance measures]\label{prop:conseuclid}
All measures consisting of the Euclidean distance of a given point $\xv$~to the cut hyperplane are dominance-consistent with respect to any set of cuts if,
for any two cuts in the set, the cut with the smallest distance measure cuts off $\xv$ and the projection of \xv~onto its hyperplane is LP-feasible.
\end{proposition}
\proof{
Let $(A)$ and $(B)$ be the two cuts of the set.
The proof directly applies to more cuts by induction.
We assume w.l.o.g.~that $d(\xv, \coefficients_A,\beta_A) > d(\xv, \coefficients_B,\beta_B)$.
The set of cuts with distance measure $d(\xv, \cdot, \cdot) = d(\xv, \coefficients_A,\beta_A)$
is a subset of the hyperplanes tangent to the sphere of radius
$d(\xv, \coefficients_A,\beta_A)$ centered at \xv.
The cut $(B)$ does not separate the projection of \xv~onto its own hyperplane, which is by assumption LP-feasible.
The projected point lies on the sphere of radius $d(\xv, \coefficients_B,\beta_B)$ centered at
\xv, and is contained in the open ball of radius $d(\xv, \coefficients_A,\beta_A)$
centered at \xv. As a tangent hyperplane to the sphere of radius
$d(\xv, \coefficients_A,\beta_A)$, $(A)$ therefore cuts off the projected point by at least $d(\xv, \coefficients_A,\beta_A) - d(\xv, \coefficients_B,\beta_B)$ and can therefore not be dominated.\qed}

Proposition~\ref{prop:conseuclid} directly applies to efficacy and analytic efficacy, with the key restriction that the projection of the relaxation point onto the cut must be LP-feasible, which excludes a majority of real instances.
Furthermore, since \noargeff is a linear function of \lpoptimal,
the property also applies to \noargavgeff as a distance measure,
with the point to project \xv being the average of the multiple LP solutions.
Finally, we can construct counter-examples where
dominance-consistency does not hold for measures based on the
Euclidean projection of a point onto the cut hyperplane in cases
where the projection is not LP-feasible, as shown in Figure~\ref{fig:infeasible_projection}.

\begin{proposition}[Consistency of the minimum efficacy]\label{prop:consmineff}
Given the set of LP solutions \lpoptimalset, we define the active solutions for a cut
$(\coefficients,\beta)$ as the subset $\argmin_{\xv\in\lpoptimalset} \eff{\coefficients}{\beta}{\xv}$.
\noargmineff is dominance-consistent with respect to
a set of cuts if for any two cuts $(A)$ and $(B)$ such that $(B)$ has a strictly lower \noargmineff,
there exists an active solution $\xv_0$ separated by $(B)$ such that its projection onto the hyperplane of $(B)$ is LP-feasible.
\end{proposition}
\proof{
Similarly to Proposition~\ref{prop:conseuclid}, the cut $(B)$ forms a tangent hyperplane to the sphere centered
at $\xv_0$ and of radius equal to the score of $(B)$.
That point is not separated by $(B)$ itself but $(A)$
has to separate it by at least the difference in score between $(A)$ and $(B)$.
\qed}

\begin{proposition}[Consistency of directed distance measures]\label{prop:consdirected}
All measures based on the distance of a point \xv~to the cut in the direction of an LP-feasible point $\hat{\xv}$~are dominance-consistent w.r.t~any set of cuts which all separate \xv.
\end{proposition}
\proof{
All points on the segment $[\xv,\hat{\xv}]$ are LP-feasible. Let $(A)$ and $(B)$ be two cuts of the set such that $(A)$ has a strictly greater measure value. Since the measure corresponds to the length of the segment cut off by cuts, some points are separated by $(A)$ only,
$(B)$ cannot dominate $(A)$.\qed
}

\begin{figure}
\centering
\begin{tikzpicture}[scale=1.0]
    \foreach \x/\y in {1/2, 1/3, 1/4,2/2, 2/3, 2/4, 3/2, 3/3, 3/4, 4/2, 4/3, 4/4}
    {
    \fill[opacity=0.4] (\x,\y) circle (2pt);
    }
    \foreach \x/\y in {2/2, 3/2, 3/3}
    {
    \fill[opacity=1] (\x,\y) circle (2pt);
    }
    \node[label=right:{\textcolor{green}{$-c$}}] (e) at (3.3, 3.5) {};
    \node[label=above right:{$\lpoptimal$}] (a) at (4, 3.2) {};
    \node (b) at (2, 2) {};
    \node (c) at (1.7778, 2.5002) {};
    \node (d) at (3, 3) {};
    \node (e) at (3, 2) {};
    
    \filldraw[thick,fill=blue!60,fill opacity=0.4] (a.center) -- (e.center) -- (b.center) -- (c.center) -- (d.center) -- cycle;
    
    \fill[draw=black,fill=red] (a.center) circle (2pt);
    
    \draw[ultra thick, dashed, red] (2.8492, 3.7538) -- (3.1987, 2.006);
    \draw[ultra thick, dotted, red] (2.2198, 3.6504) -- (3.978, 2.1855);
    \draw[ultra thick, green] (3.5, 3.2) -- (4.5, 3.2);
\end{tikzpicture}
\caption{Two cuts with LP-feasible orthogonal projections of \lpoptimal. The dotted cut is dominated, but has a better \noargexpimprov score.}
\label{fig:expimprov-nonconsist}
\end{figure}
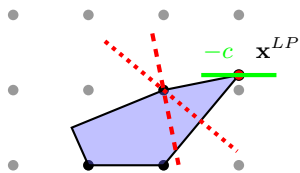

Proposition~\ref{prop:consdirected} notably applies to directed cutoff distance and
(approximate) analytic directed cutoff distance. We note that dominance-consistency does not extend to \noargexpimprov, even with feasible LP projections, with a counterexample visualised in Figure \ref{fig:expimprov-nonconsist}.

\subsection{Multi-Output Regression}
\label{subsec:ml_intro}

Machine Learning for MIP has mainly focused on classification tasks, e.g.~should
an algorithm be run with option A or B, or run at all.
In our case, the output of interest is the (relative) advantage of distance measures in terms of some performance criterion.
For some instances, different measures may result in (near-)identical performance,
e.g.~if the same, ``obvious'' subset of cuts is selected by all methods.
Classifying these (near-)ties with one distance measure or choosing an arbitrary threshold to classify a measure as well-performing could prevent
the model from fitting well on the important data points where the selection method is significant.
We therefore pose our learning task as a multi-output regression that predicts the relative performance of each method.
We aim at an interpretable model to predict the preferred measure, with the goal of outperforming individual distance measures.

The feature space used as input to our model consists of: \emph{dual degeneracy} (fraction of non-basic variables with zero-reduced cost), \emph{primal degeneracy} (fraction of basic variables at their respective bounds), \emph{solution fractionality} (fraction of integer variables with fractional LP values), \emph{thinness} (fraction of equality constraints), and \emph{density} (fraction of non-zero entries in constraint matrix). All features are obtainable at the root node before the separation process begins, are relevant to the separation process, and are easy to retrieve.

\section{Experiments}

We perform experiments on the MIPLIB 2017 collection set \footnote{MIPLIB 2017 -- The Mixed Integer Programming Library \url{https://miplib.zib.de/}.} \cite{miplib}, which we will now simply refer to as MIPLIB. We define a run as an instance random-seed pair for which we enforce exactly 50 separation rounds at the root node, with a maximum of 10 cuts to be added per round. We use default separators, but increase the amount of cuts that can be generated. Additionally, in order to reduce the variability of the solving process, restarts are disabled, no cuts are allowed to be added after the root node, and the best available MIPLIB solution is provided. All other aspects of the solver are untouched, with our experiments only replacing the cut scoring function in SCIP, not the selection algorithm itself, see \cite{turner2022adaptive}. The primary motivation is to assess the performance of different distance-based cut measures, and to determine for which instance characteristics a distance measure is effective. All results are obtained by averaging results over the SCIP random seeds $\{1,2,3\}$, and instances are filtered subject to Table~\ref{tab:instance_filtering}, with 162 instances remaining. Three LP solutions are used for \noargmineff and \noargavgeff.

\begin{table}[h]
\scriptsize
\centering
\resizebox{\columnwidth}{!}{%
\begin{tabular}{lc}
Criteria & \% of instances removed \\
\hline
Tags: \textit{feasibility}, \textit{numerics}, \textit{infeasible}, \textit{no solution} & 4.5\%, 17.5\%, 2.8\%, 0.9\% \\
Unbounded objective, MIPLIB solution unavailable & 0.9\%, 2.6\% \\
Root optimal (any measure) & 13.3\% \\
Root node with separation rounds longer than 600s (any measure) & 13.4\% \\
No optional cuts generated (all measures) & 2.3\% \\
Numerical issues (any measure) & 1.2\% \\
Failed to prove optimality in branch and cut within 7200s (all measures) & 21.7\% \\
\hline
\end{tabular}
}
\caption{Criteria for which we removed instances from the MIPLIB collection and percentages of instances affected by each criterion}
\label{tab:instance_filtering}
\end{table}

For all experiments, SCIP 8.0.2 \cite{scip8} is used, with PySCIPOpt \cite{pyscipopt} as the API, and Xpress 8.14 \cite{xpress} as the LP solver. All experiments are run on a cluster equipped with Intel Xeon Gold 5122 CPUs with 3.60GHz and 96GB main memory. The code used for all experiments is available and open-source\footnote{\url{https://github.com/Opt-Mucca/Analytic-Center-Cut-Selection}}.
The structure of this section is as follows. In Subsection~\ref{subsec:root_experiments}, we present results of our distance-based cut measures on root node restricted runs. In Subsection~\ref{subsec:tree_experiments} we present results of our distance measures generalised to branch and cut. Finally, in Subsection~\ref{subsec:ml_experiments} we present the performance of our support vector regression model on selecting distance measures.

\subsection{Root Node Results}
\label{subsec:root_experiments}

\begin{table}[h]
\small
\centering
\begin{tabular}{lll}
Function & Measure & Description \\
\hline
\noargeff & Efficacy & See Section \ref{sec:introduction} \\
\noargdcd & Directed cutoff distance & See Section \ref{sec:introduction} \\
\noarganaleff & Analytic Efficacy & See Subsection \ref{subsec:ae} \\
\noarganaldcd & Analytic directed cutoff distance & See Subsection \ref{subsec:ae} \\
\noargapproxanaldcd & Approximate analytic directed cutoff distance & See Subsection \ref{subsec:ae} \\
\noargavgeff & Average efficacy & See Subsection \ref{subsec:multiple_lp} \\
\noargmineff & Minimum efficacy & See Subsection \ref{subsec:multiple_lp} \\
\noargexpimprov & Expected Improvement & See Section \ref{sec:introduction} \\
\hline
\end{tabular}
\caption{Summary of all distance measures}
\label{tab:cutoff_methods}
\end{table}

Table~\ref{tab:cutoff_methods} provides a summary of all cut selection measures we evaluated.
We compare head-to-head results on the primal-dual difference after 50 separation rounds. We say that a scoring measure has outperformed another for an instance, if it is at least as good over all random seeds, and strictly better for at least one. Curiously, we observe a clear hierarchy of distance-based cut measures for root-restricted dual bound performance over MIPLIB. That is, $\noarganaldcd \geq \noargapproxanaldcd \geq \noarganaleff \geq \noargmineff \geq \noargavgeff \geq \noargdcd \geq \noargeff \geq \noargexpimprov$. 

\begin{table}[]
    \scriptsize
    \centering
    \resizebox{\columnwidth}{!}{%
    \begin{tabular}{l|c|c|c|c|c|c|c|c}
    & \noarganaldcd & \noargapproxanaldcd & \noarganaleff & \noargmineff & \noargavgeff & \noargeff & \noargdcd & \noargexpimprov \\
    \noarganaldcd & - & \textbf{\textcolor{blue}{0.29}/\textcolor{blue}{0.16}} & \textbf{\textcolor{blue}{0.32}/\textcolor{blue}{0.3}} & 0.29/0.29 & \textbf{\textcolor{blue}{0.32}/\textcolor{blue}{0.26}} & \textbf{\textcolor{blue}{0.38}/\textcolor{blue}{0.22}} & \textbf{\textcolor{blue}{0.4}/\textcolor{blue}{0.21}} & \textbf{\textcolor{blue}{0.47}/\textcolor{blue}{0.22}} \\
    \noargapproxanaldcd & \textsl{\textcolor{red}{0.16}/\textcolor{red}{0.29}} & - & \textbf{\textcolor{blue}{0.3}/\textcolor{blue}{0.3}} & \textbf{\textcolor{blue}{0.3}/\textcolor{blue}{0.27}} & \textbf{\textcolor{blue}{0.3}/\textcolor{blue}{0.27}} & \textbf{\textcolor{blue}{0.38}/\textcolor{blue}{0.22}} & \textbf{\textcolor{blue}{0.37}/\textcolor{blue}{0.23}} & \textbf{\textcolor{blue}{0.43}/\textcolor{blue}{0.25}} \\
    \noarganaleff & \textsl{\textcolor{red}{0.3}/\textcolor{red}{0.32}} & \textsl{\textcolor{red}{0.3}/\textcolor{red}{0.3}} & - & \textbf{\textcolor{blue}{0.25}/\textcolor{blue}{0.23}} & \textbf{\textcolor{blue}{0.22}/\textcolor{blue}{0.2}} & \textbf{\textcolor{blue}{0.31}/\textcolor{blue}{0.14}} & \textbf{\textcolor{blue}{0.37}/\textcolor{blue}{0.26}} & \textbf{\textcolor{blue}{0.41}/\textcolor{blue}{0.23}} \\
    \noargmineff & 0.29/0.29 & \textsl{\textcolor{red}{0.27}/\textcolor{red}{0.3}} & \textsl{\textcolor{red}{0.23}/\textcolor{red}{0.25}} & - & \textbf{\textcolor{blue}{0.14}/\textcolor{blue}{0.13}} & \textbf{\textcolor{blue}{0.28}/\textcolor{blue}{0.12}} & \textbf{\textcolor{blue}{0.32}/\textcolor{blue}{0.26}} & \textbf{\textcolor{blue}{0.36}/\textcolor{blue}{0.24}} \\
    \noargavgeff & \textsl{\textcolor{red}{0.26}/\textcolor{red}{0.32}} & \textsl{\textcolor{red}{0.27}/\textcolor{red}{0.3}} & \textsl{\textcolor{red}{0.2}/\textcolor{red}{0.22}} & \textsl{\textcolor{red}{0.13}/\textcolor{red}{0.14}} & - & \textbf{\textcolor{blue}{0.27}/\textcolor{blue}{0.12}} & \textbf{\textcolor{blue}{0.31}/\textcolor{blue}{0.25}} & \textbf{\textcolor{blue}{0.38}/\textcolor{blue}{0.24}} \\
    \noargeff & \textsl{\textcolor{red}{0.22}/\textcolor{red}{0.38}} & \textsl{\textcolor{red}{0.22}/\textcolor{red}{0.38}} & \textsl{\textcolor{red}{0.14}/\textcolor{red}{0.31}} & \textsl{\textcolor{red}{0.12}/\textcolor{red}{0.28}} & \textsl{\textcolor{red}{0.12}/\textcolor{red}{0.27}} & - & \textbf{\textcolor{blue}{0.26}/\textcolor{blue}{0.25}} & \textbf{\textcolor{blue}{0.32}/\textcolor{blue}{0.25}} \\
    \noargdcd & \textsl{\textcolor{red}{0.21}/\textcolor{red}{0.4}} & \textsl{\textcolor{red}{0.23}/\textcolor{red}{0.37}} & \textsl{\textcolor{red}{0.26}/\textcolor{red}{0.37}} & \textsl{\textcolor{red}{0.26}/\textcolor{red}{0.32}} & \textsl{\textcolor{red}{0.25}/\textcolor{red}{0.31}} & \textsl{\textcolor{red}{0.25}/\textcolor{red}{0.26}} & - & \textbf{\textcolor{blue}{0.36}/\textcolor{blue}{0.26}} \\
    \noargexpimprov & \textsl{\textcolor{red}{0.22}/\textcolor{red}{0.47}} & \textsl{\textcolor{red}{0.25}/\textcolor{red}{0.43}} & \textsl{\textcolor{red}{0.23}/\textcolor{red}{0.41}} & \textsl{\textcolor{red}{0.24}/\textcolor{red}{0.36}} & \textsl{\textcolor{red}{0.24}/\textcolor{red}{0.38}} & \textsl{\textcolor{red}{0.25}/\textcolor{red}{0.32}} & \textsl{\textcolor{red}{0.26}/\textcolor{red}{0.36}} & - 
    \end{tabular}
    }
    \caption{Entry coordinate $(i,j)$ is a tuple of win / loss percentage over all instances for dual bound improvement of measure $i$ over measure $j$. A win is defined by at least as good performance over all seeds, and better performance for at least one seed.}
    \label{tab:pairwise_root_wins}
\end{table}

\ifarxiv
\else
\vspace{-1.2cm}
\fi

\begin{figure}[]
    \centering
    \includegraphics[scale=0.4]{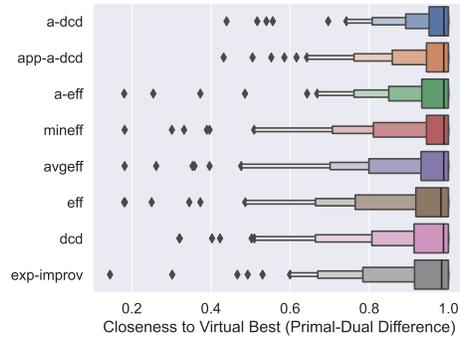}
    \caption{Boxenplots of distance measures root-node performance. Each instance compared to the virtual best and averaged over random seeds.}
    \label{fig:root_primal_dual_boxenplots}
\end{figure}

A head-to-head comparison, while helpful for ranking measures, contains limited information on the distribution of performance. For this reason, we also visualise the primal-dual difference results using boxenplots in Figure \ref{fig:root_primal_dual_boxenplots}, see \cite{letter-value-plot} for a description. Figure \ref{fig:root_primal_dual_boxenplots} shows a comparison of each method against the so-called \emph{virtual best}. Therefore, we divide, instance by instance, the average gap (over all seeds) by the best average gap among the eight methods. We observe similar results to Table \ref{tab:pairwise_root_wins} in that \noarganaldcd, \noargapproxanaldcd, and \noarganaleff are superior to other methods in terms of dual bound improvement. It should be noted, however, that on average \noarganaleff takes 32\% of the root node processing time and \noarganaldcd takes 25\%. This is in contrast to \noargeff, which only takes 0.8\%. We conclude that using an analytic center for cut selection is beneficial for closing the primal-dual gap during root node cutting. 



\subsection{Branch and Bound Generalisation}
\label{subsec:tree_experiments}

\begin{table}[h]
    \scriptsize
    \begin{subtable}{1.0\textwidth}
    \centering
    \begin{tabular}{l|c|c|c|c|c|c|c|c}
    & \noarganaldcd & \noargapproxanaldcd & \noarganaleff & \noargmineff & \noargavgeff & \noargeff & \noargdcd & \noargexpimprov \\
    \noarganaldcd & - & \textbf{\textcolor{blue}{0.17}/\textcolor{blue}{0.16}} & \textbf{\textcolor{blue}{0.29}/\textcolor{blue}{0.17}} & \textbf{\textcolor{blue}{0.25}/\textcolor{blue}{0.17}} & \textbf{\textcolor{blue}{0.23}/\textcolor{blue}{0.19}} & \textbf{\textcolor{blue}{0.27}/\textcolor{blue}{0.17}} & \textbf{\textcolor{blue}{0.23}/\textcolor{blue}{0.17}} & \textbf{\textcolor{blue}{0.28}/\textcolor{blue}{0.15}} \\
    \noargapproxanaldcd & \textsl{\textcolor{red}{0.16}/\textcolor{red}{0.17}} & - & \textsl{\textcolor{red}{0.2}/\textcolor{red}{0.22}} & \textbf{\textcolor{blue}{0.23}/\textcolor{blue}{0.2}} & \textsl{\textcolor{red}{0.19}/\textcolor{red}{0.27}} & \textbf{\textcolor{blue}{0.22}/\textcolor{blue}{0.19}} & \textsl{\textcolor{red}{0.21}/\textcolor{red}{0.22}} & \textbf{\textcolor{blue}{0.23}/\textcolor{blue}{0.18}} \\
    \noarganaleff & \textsl{\textcolor{red}{0.17}/\textcolor{red}{0.29}} & \textbf{\textcolor{blue}{0.22}/\textcolor{blue}{0.2}} & - & 0.15/0.15 & 0.2/0.2 & 0.17/0.17 & \textbf{\textcolor{blue}{0.2}/\textcolor{blue}{0.19}} & \textbf{\textcolor{blue}{0.18}/\textcolor{blue}{0.17}} \\
    \noargmineff & \textsl{\textcolor{red}{0.17}/\textcolor{red}{0.25}} & \textsl{\textcolor{red}{0.2}/\textcolor{red}{0.23}} & 0.15/0.15 & - & \textsl{\textcolor{red}{0.11}/\textcolor{red}{0.17}} & \textbf{\textcolor{blue}{0.17}/\textcolor{blue}{0.1}} & \textsl{\textcolor{red}{0.17}/\textcolor{red}{0.18}} & \textbf{\textcolor{blue}{0.18}/\textcolor{blue}{0.15}} \\
    \noargavgeff & \textsl{\textcolor{red}{0.19}/\textcolor{red}{0.23}} & \textbf{\textcolor{blue}{0.27}/\textcolor{blue}{0.19}} & 0.2/0.2 & \textbf{\textcolor{blue}{0.17}/\textcolor{blue}{0.11}} & - & \textbf{\textcolor{blue}{0.17}/\textcolor{blue}{0.14}} & \textbf{\textcolor{blue}{0.21}/\textcolor{blue}{0.19}} & \textbf{\textcolor{blue}{0.23}/\textcolor{blue}{0.22}} \\
    \noargeff & \textsl{\textcolor{red}{0.17}/\textcolor{red}{0.27}} & \textsl{\textcolor{red}{0.19}/\textcolor{red}{0.22}} & 0.17/0.17 & \textsl{\textcolor{red}{0.1}/\textcolor{red}{0.17}} & \textsl{\textcolor{red}{0.14}/\textcolor{red}{0.17}} & - & \textsl{\textcolor{red}{0.16}/\textcolor{red}{0.2}} & \textsl{\textcolor{red}{0.19}/\textcolor{red}{0.2}} \\
    \noargdcd & \textsl{\textcolor{red}{0.17}/\textcolor{red}{0.23}} & \textbf{\textcolor{blue}{0.22}/\textcolor{blue}{0.21}} & \textsl{\textcolor{red}{0.19}/\textcolor{red}{0.2}} & \textbf{\textcolor{blue}{0.18}/\textcolor{blue}{0.17}} & \textsl{\textcolor{red}{0.19}/\textcolor{red}{0.21}} & \textbf{\textcolor{blue}{0.2}/\textcolor{blue}{0.16}} & - & \textbf{\textcolor{blue}{0.19}/\textcolor{blue}{0.17}} \\
    \noargexpimprov & \textsl{\textcolor{red}{0.15}/\textcolor{red}{0.28}} & \textsl{\textcolor{red}{0.18}/\textcolor{red}{0.23}} & \textsl{\textcolor{red}{0.17}/\textcolor{red}{0.18}} & \textsl{\textcolor{red}{0.15}/\textcolor{red}{0.18}} & \textsl{\textcolor{red}{0.22}/\textcolor{red}{0.23}} & \textbf{\textcolor{blue}{0.2}/\textcolor{blue}{0.19}} & \textsl{\textcolor{red}{0.17}/\textcolor{red}{0.19}} & - 
    \end{tabular}
    \caption{Number of nodes}
    \label{tab:pairwise_node_wins}
    \end{subtable}
    \begin{subtable}{1.0\textwidth}
    \centering
    \begin{tabular}{l|c|c|c|c|c|c|c|c}
    & \noarganaldcd & \noargapproxanaldcd & \noarganaleff & \noargmineff & \noargavgeff & \noargeff & \noargdcd & \noargexpimprov \\
    \noarganaldcd & - & \textbf{\textcolor{blue}{0.37}/\textcolor{blue}{0.14}} & \textbf{\textcolor{blue}{0.33}/\textcolor{blue}{0.16}} & \textbf{\textcolor{blue}{0.28}/\textcolor{blue}{0.15}} & \textbf{\textcolor{blue}{0.28}/\textcolor{blue}{0.17}} & \textsl{\textcolor{red}{0.22}/\textcolor{red}{0.3}} & \textsl{\textcolor{red}{0.16}/\textcolor{red}{0.3}} & \textsl{\textcolor{red}{0.19}/\textcolor{red}{0.29}} \\
    \noargapproxanaldcd & \textsl{\textcolor{red}{0.14}/\textcolor{red}{0.37}} & - & \textbf{\textcolor{blue}{0.23}/\textcolor{blue}{0.2}} & \textbf{\textcolor{blue}{0.25}/\textcolor{blue}{0.23}} & \textbf{\textcolor{blue}{0.24}/\textcolor{blue}{0.23}} & \textsl{\textcolor{red}{0.17}/\textcolor{red}{0.35}} & \textsl{\textcolor{red}{0.15}/\textcolor{red}{0.3}} & \textsl{\textcolor{red}{0.16}/\textcolor{red}{0.36}} \\
    \noarganaleff & \textsl{\textcolor{red}{0.16}/\textcolor{red}{0.33}} & \textsl{\textcolor{red}{0.2}/\textcolor{red}{0.23}} & - & \textbf{\textcolor{blue}{0.23}/\textcolor{blue}{0.23}} & \textsl{\textcolor{red}{0.26}/\textcolor{red}{0.29}} & \textsl{\textcolor{red}{0.14}/\textcolor{red}{0.41}} & \textsl{\textcolor{red}{0.15}/\textcolor{red}{0.35}} & \textsl{\textcolor{red}{0.1}/\textcolor{red}{0.37}} \\
    \noargmineff & \textsl{\textcolor{red}{0.15}/\textcolor{red}{0.28}} & \textsl{\textcolor{red}{0.23}/\textcolor{red}{0.25}} & \textsl{\textcolor{red}{0.23}/\textcolor{red}{0.23}} & - & 0.2/0.2 & \textsl{\textcolor{red}{0.12}/\textcolor{red}{0.4}} & \textsl{\textcolor{red}{0.11}/\textcolor{red}{0.33}} & \textsl{\textcolor{red}{0.12}/\textcolor{red}{0.29}} \\
    \noargavgeff & \textsl{\textcolor{red}{0.17}/\textcolor{red}{0.28}} & \textsl{\textcolor{red}{0.23}/\textcolor{red}{0.24}} & \textbf{\textcolor{blue}{0.29}/\textcolor{blue}{0.26}} & 0.2/0.2 & - & \textsl{\textcolor{red}{0.1}/\textcolor{red}{0.37}} & \textsl{\textcolor{red}{0.12}/\textcolor{red}{0.34}} & \textsl{\textcolor{red}{0.12}/\textcolor{red}{0.32}} \\
    \noargeff & \textbf{\textcolor{blue}{0.3}/\textcolor{blue}{0.22}} & \textbf{\textcolor{blue}{0.35}/\textcolor{blue}{0.17}} & \textbf{\textcolor{blue}{0.41}/\textcolor{blue}{0.14}} & \textbf{\textcolor{blue}{0.4}/\textcolor{blue}{0.12}} & \textbf{\textcolor{blue}{0.37}/\textcolor{blue}{0.1}} & - & \textsl{\textcolor{red}{0.18}/\textcolor{red}{0.22}} & 0.2/0.2 \\
    \noargdcd & \textbf{\textcolor{blue}{0.3}/\textcolor{blue}{0.16}} & \textbf{\textcolor{blue}{0.3}/\textcolor{blue}{0.15}} & \textbf{\textcolor{blue}{0.35}/\textcolor{blue}{0.15}} & \textbf{\textcolor{blue}{0.33}/\textcolor{blue}{0.11}} & \textbf{\textcolor{blue}{0.34}/\textcolor{blue}{0.12}} & \textbf{\textcolor{blue}{0.22}/\textcolor{blue}{0.18}} & - & \textbf{\textcolor{blue}{0.21}/\textcolor{blue}{0.2}} \\
    \noargexpimprov & \textbf{\textcolor{blue}{0.29}/\textcolor{blue}{0.19}} & \textbf{\textcolor{blue}{0.36}/\textcolor{blue}{0.16}} & \textbf{\textcolor{blue}{0.37}/\textcolor{blue}{0.1}} & \textbf{\textcolor{blue}{0.29}/\textcolor{blue}{0.12}} & \textbf{\textcolor{blue}{0.32}/\textcolor{blue}{0.12}} & 0.2/0.2 & \textsl{\textcolor{red}{0.2}/\textcolor{red}{0.21}} & -
    \end{tabular}
    \caption{Time}
    \label{tab:pairwise_time_wins}
    \end{subtable}
    \caption{Entry coordinate $(i,j)$ is a tuple of win / loss percentage over all instances for measure $i$ over measure $j$. A win is defined by at least as good performance over all seeds, and better performance for at least one seed.}
    \label{tab:headtoheadtree}
\end{table}
\ifarxiv\else\vspace*{-1cm}\fi
While root node performance can be used as a surrogate for solver performance, there is no guarantee that results generalise to the entire solving process. We therefore extend our experiments to the branch-and-bound tree with a time limit of two hours. The head-to-head results of each distance measure for the number of nodes and solve time are displayed in Table~\ref{tab:headtoheadtree}.
For the number of nodes, note that we removed all instances where a measure timed out. We observe in the number of nodes comparison that \noarganaldcd remains the superior method, however the ordering of methods is now less clear. Most interesting is the drop in performance of \noargapproxanaldcd compared to the root node results, which suggests that the analytic center from previous separation rounds is often not a good direction for distance measures. This is supported by the fact that $26.5\%$ of the analytic centers from previous rounds are LP infeasible. In the solve time comparison, we see that the `cheaper' measures, \noargeff, \noargdcd, \noargexpimprov, which require no additional LP solver calls, are superior over the more `expensive' measures. This suggests that while our introduced methods, especially \noarganaldcd, can reduce the number of nodes and have better root node performance, the total solve time is not similarly improved. We note that while \noargdcd is superior in Table~\ref{tab:pairwise_time_wins}, we believe that our experimental design is overly favourable since we start with an optimal MIPLIB solution (or best known, for unsolved instances).

\begin{figure}[htp]
  \centering
  \includegraphics[height=5cm,width=6.3cm]{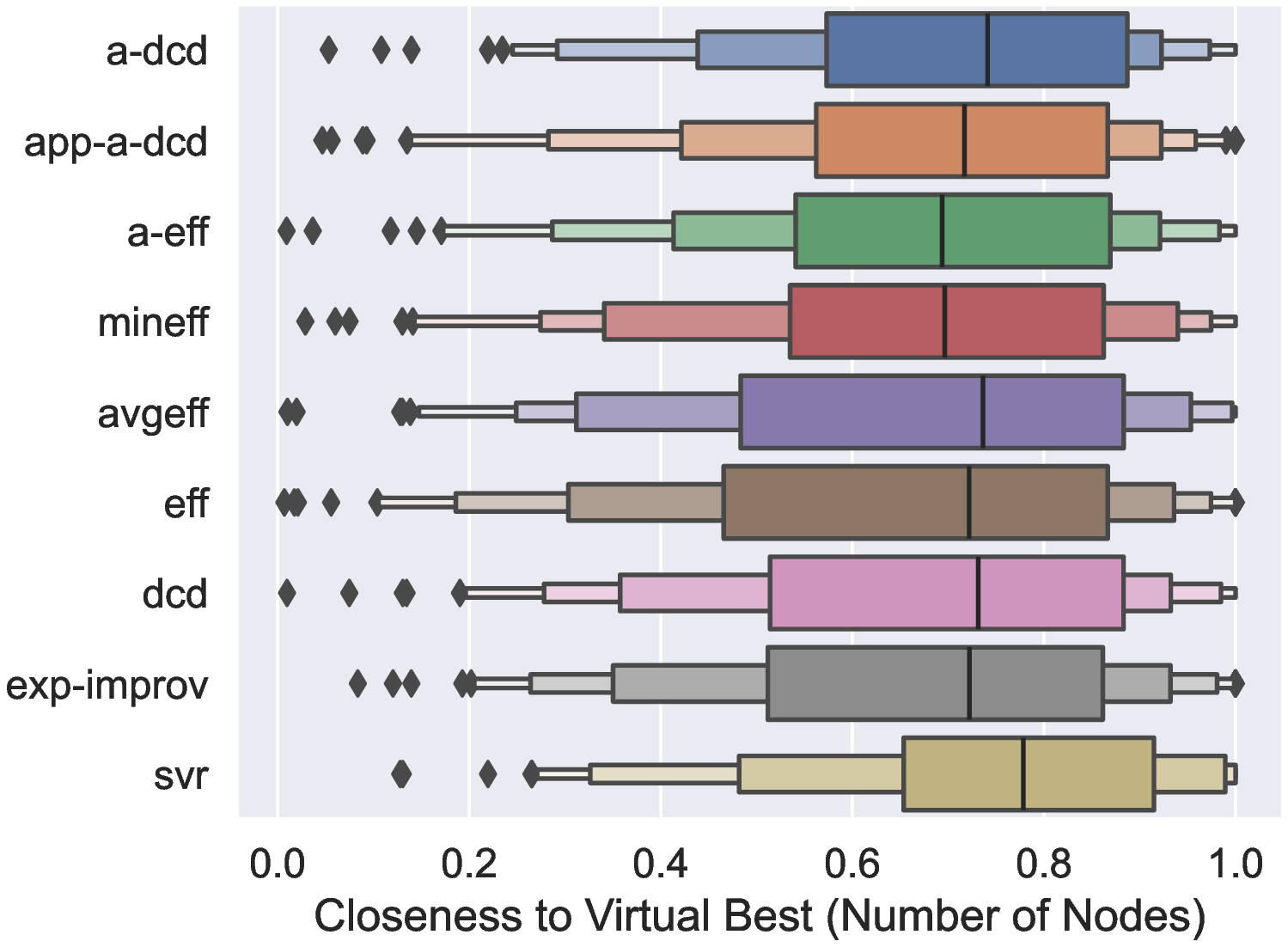}\hspace{.1em}%
  \includegraphics[height=5cm,width=5.6cm]{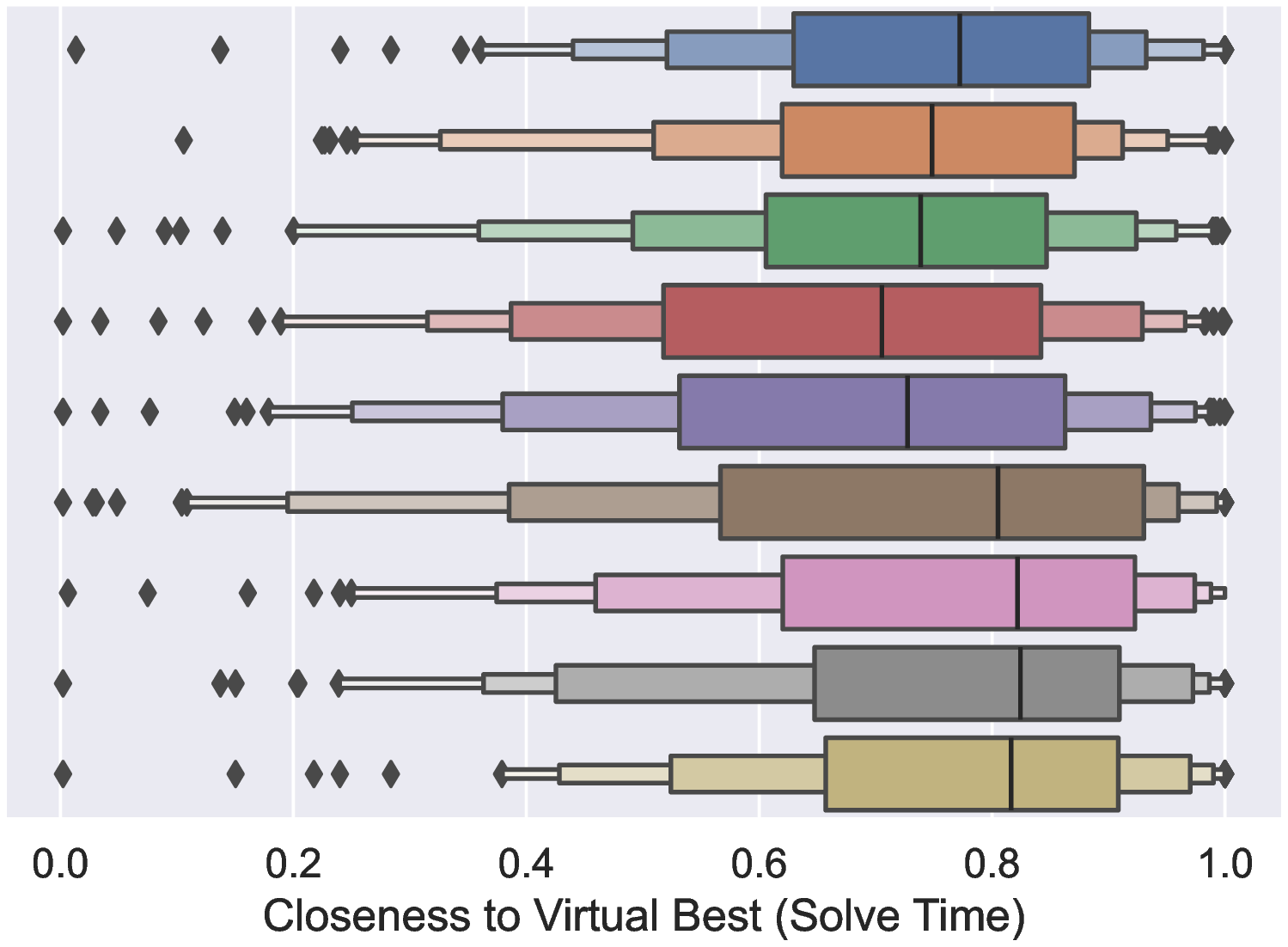}
  \caption{Boxenplots of measures' tree performance. Nodes (left), time (right).}
  \label{fig:tree_node_time_boxenplot}
\end{figure}

Similarly to Subsection~\ref{subsec:root_experiments}, we visualise an instance-wide comparison to the virtual best of all measures for number of nodes and solve time in Figure~\ref{fig:tree_node_time_boxenplot}.
The results confirm our conclusion from Table~\ref{tab:pairwise_node_wins} that \noarganaldcd is the best performing measure, and that \noargeff is the worst one w.r.t.~number of nodes. All other measures, however, have similar distributions, making stronger conclusions difficult. We note that $90.4\%$ of cuts have infeasible projections when scoring by \noargeff, and that previous studies, see \cite{mipdegeneracy}, identify an 87.5\% occurrence rate of some level of dual degeneracy of the final root node LP in standard benchmark instances, confirming the practical geometric limitations of efficacy presented in Section~\ref{subsec:properties}. For solve time, we observe the improved performance of `cheaper' methods through their relatively high median values. We also observe that \noarganaldcd has the smallest performance variability of all measures, while the standard \noargeff has the largest performance variability. This implies that using an analytic center for cut selection can help to reduce performance variability, an interesting observation by its own.


\subsection{Regression Model Results}
\label{subsec:ml_experiments}

We have thus far observed that our newly introduced measures, especially \noarganaldcd, have superior root node dual bound performance than traditional measures, and often result in smaller branch-and-bound trees. No single measure is however dominant, as seen in Tables~\ref{tab:pairwise_node_wins} and \ref{tab:pairwise_time_wins}, with no single measure ever having less than 10\% of instances as wins in the head-to-head contest. This motivates the need for an adaptive method, which decides on a distance measure at the start of the solving process that will best perform on the instance.

We use support vector regression (SVR) with a cubic kernel function, see \cite{smola2004tutorial},
implemented in \texttt{scikit-learn}~\cite{scikit-learn} with default parameters. We train on instance-seed pairs, with the virtual best number of nodes for each pair divided by the number of nodes under the distance measure as output. Our model was trained using 5-fold cross-validation, with 10\% of pairs retained for validation.
We were able to achieve comparable performance with regression forests and alternative kernels, likewise with default parameter sets. The final model was selected due to its ease of interpretation and potential embedding in a MIP solver.

\begin{figure}[h]
    \centering
    \includegraphics[width=0.73\linewidth]{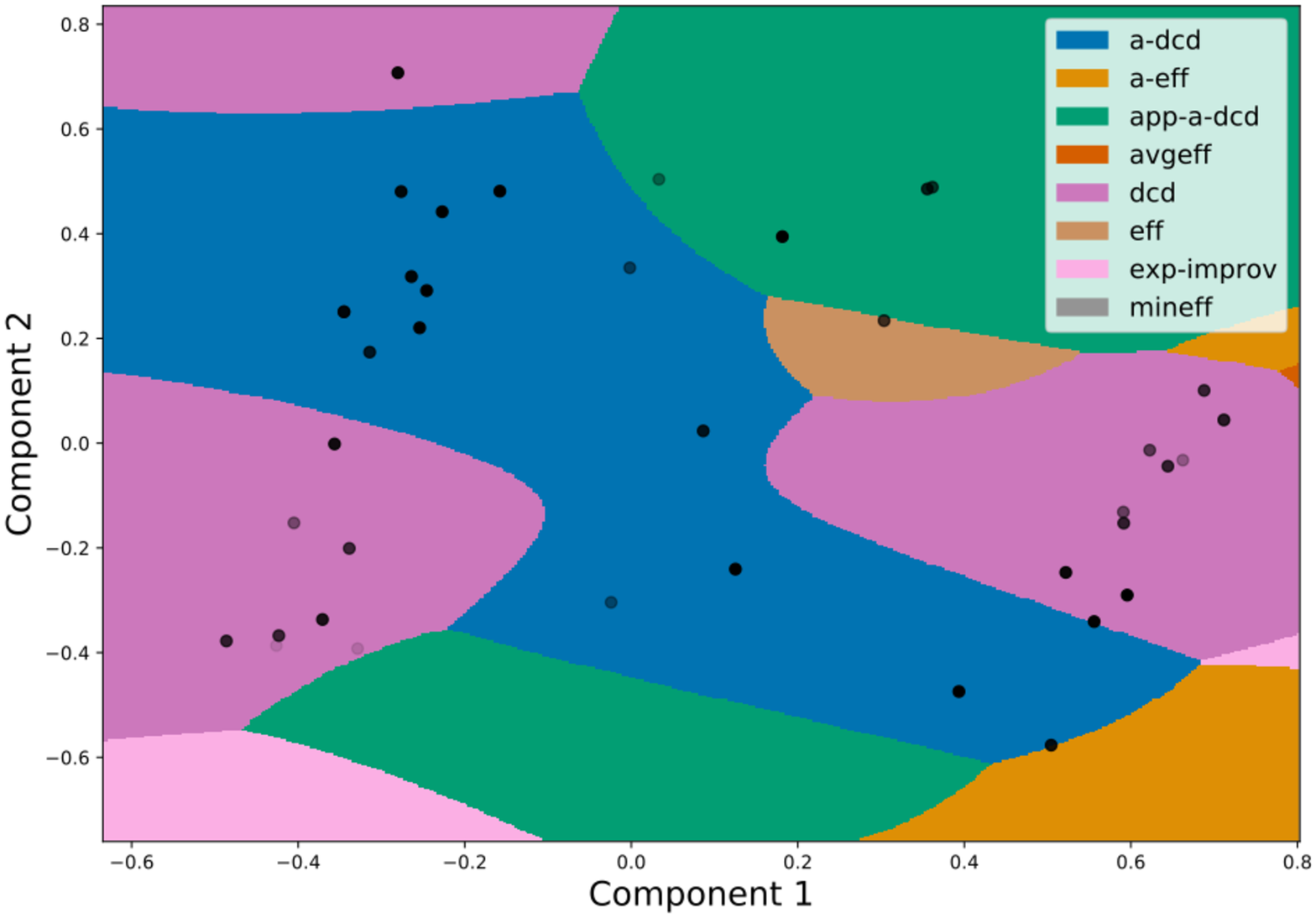}
    \caption[.]{Decision regions in transformed feature space. Dots are validation instances, with their opacity the relative performance of the predicted measure.\\
    {\scriptsize
    $\text{Component 1:}\,\,\, 0.947 \text{ dual\_deg} - 0.205 \text{ primal\_deg} + 0.22 \text{ frac} - 0.089 \text{ thin} - 0.063 \text{ density}$\\
    $\text{Component 2:}\,\, -0.27 \text{ dual\_deg} - 0.733 \text{ primal\_deg} + 0.467 \text{ frac} - 0.256 \text{ thin} + 0.326 \text{ density}$
    }
    }
    \label{fig:pca_space}
\end{figure}

We observe in Figure~\ref{fig:tree_node_time_boxenplot} that our trained model clearly outperforms any individual distance measure w.r.t.~number of nodes. Further, when considering the shifted geometric mean of the number of nodes, it is 12\% smaller than that of the best overall performing distance measure.
This strong result does not generalise to solve time, however. The distribution looks similar to \noarganaldcd, albeit with a better median. The shifted geometric mean of our model w.r.t.~solve time is 8\% larger than the single best-performing distance measure. Note that there are situations where available memory is a limiting factor -- e.g., super-computing -- which makes node savings important.
Finally, we visualise the decision boundaries of the trained model over the two first principal components from a PCA of the original features, maintaining 71\% of the explained variance.
We determine decision boundaries with the largest regression value over all distances and visualise the result in Figure~\ref{fig:pca_space}, with the component equations printed below.
\ifarxiv
\else
\vspace{-1cm}
\fi

\ifarxiv\else\vspace*{-1cm}\fi
\section{Conclusion}
In this paper, we reassessed the question of cut selection through the lens of distance measures. Motivated by geometric properties of polyhedra encountered in MIPs, we defined measures based on analytic centers and multiple LP solutions. We showed their performance, and more importantly, that the relative performance of distance measures can be learned for new instances with an interpretable and implementable model. We found that the introduced measures help to reduce root node gap, size of the branch and bound tree and performance variability.
The focus of our work is on the improved evaluation of individual cuts; future directions will build upon these measures to enhance the whole separation process, incorporating combinations of cut measures, and generalising the promising node reductions to improved runtime.

\section*{Acknowledgements}
The work for this article has been conducted in the Research Campus MODAL funded by the German Federal Ministry
of Education and Research (BMBF) (fund numbers 05M14ZAM, 05M20ZBM). The described research activities are
funded by the Federal Ministry for Economic Affairs and Energy within the project UNSEEN (ID: 03EI1004-C).


\bibliographystyle{splncs04}
\bibliography{mybib}

\begin{thebibliography}{10}
\providecommand{\url}[1]{\texttt{#1}}
\providecommand{\urlprefix}{URL }
\providecommand{\doi}[1]{https://doi.org/#1}

\bibitem{achterberg2007constraint}
Achterberg, T.: Constraint integer programming. Ph.D. thesis, TU Berlin (2007)

\bibitem{cplexpatent}
Achterberg, T.: {LP} relaxation modification and cut selection in a {MIP}
  solver. {US Patent US8463729B2} (2013),
  \url{https://patents.google.com/patent/US8463729B2/en}

\bibitem{achterberg2013mixed}
Achterberg, T., Wunderling, R.: Mixed integer programming: Analyzing 12 years
  of progress. In: Facets of combinatorial optimization, pp. 449--481. Springer
  (2013)

\bibitem{zerohalf}
Andreello, G., Caprara, A., Fischetti, M.: Embedding $\{$0, $1/2\}$-cuts in a
  branch-and-cut framework: A computational study. INFORMS Journal on Computing
   \textbf{19}(2),  229--238 (2007)

\bibitem{BAENA2011310}
Baena, D., Castro, J.: Using the analytic center in the feasibility pump.
  Operations Research Letters  \textbf{39}(5),  310--317 (2011).
  \doi{10.1016/j.orl.2011.07.005},
  \url{https://www.sciencedirect.com/science/article/pii/S0167637711000824}

\bibitem{balcan2021sample}
Balcan, M.F.F., Prasad, S., Sandholm, T., Vitercik, E.: Sample complexity of
  tree search configuration: Cutting planes and beyond. Advances in Neural
  Information Processing Systems  \textbf{34} (2021)

\bibitem{baltean2019scoring}
Baltean-Lugojan, R., Bonami, P., Misener, R., Tramontani, A.: Scoring positive
  semidefinite cutting planes for quadratic optimization via trained neural
  networks. optimization-online preprint 2018/11/6943  (2019)

\bibitem{berthold2018four}
Berthold, T., Perregaard, M., M{\'e}sz{\'a}ros, C.: {Four good reasons to use
  an interior point solver within a MIP solver}. In: Operations Research
  Proceedings 2017, pp. 159--164. Springer (2018)

\bibitem{scip8}
Bestuzheva, K., Besançon, M., Chen, W.K., Chmiela, A., Donkiewicz, T., van
  Doornmalen, J., Eifler, L., Gaul, O., Gamrath, G., Gleixner, A., Gottwald,
  L., Graczyk, C., Halbig, K., Hoen, A., Hojny, C., van~der Hulst, R., Koch,
  T., Lübbecke, M., Maher, S.J., Matter, F., Mühmer, E., Müller, B.,
  Pfetsch, M.E., Rehfeldt, D., Schlein, S., Schlösser, F., Serrano, F.,
  Shinano, Y., Sofranac, B., Turner, M., Vigerske, S., Wegscheider, F.,
  Wellner, P., Weninger, D., Witzig, J.: {The SCIP Optimization Suite 8.0}
  (2021)

\bibitem{chmiela2022implementation}
Chmiela, A., Mu{\~n}oz, G., Serrano, F.: On the implementation and
  strengthening of intersection cuts for qcqps. Mathematical Programming pp.
  1--38 (2022)

\bibitem{dey2018theoretical}
Dey, S.S., Molinaro, M.: Theoretical challenges towards cutting-plane
  selection. Mathematical Programming  \textbf{170}(1),  237--266 (2018)

\bibitem{xpress}
{FICO Xpress Optimization}.
  \url{https://www.fico.com/en/products/fico-xpress-optimization}, accessed:
  2022-11-10

\bibitem{multiplelpcutting}
Fischetti, M., Lodi, A., Monaci, M., Salvagnin, D., Tramontani, A.: Improving
  branch-and-cut performance by random sampling. Mathematical Programming
  Computation  \textbf{8}(1),  113--132 (2016)

\bibitem{fischetti2009yoyo}
Fischetti, M., Salvagnin, D.: Yoyo search: a bisection cutting-plane method
  (2009)

\bibitem{mipdegeneracy}
Gamrath, G., Berthold, T., Salvagnin, D.: An exploratory computational analysis
  of dual degeneracy in mixed-integer programming. EURO Journal on
  Computational Optimization  \textbf{8}(3-4),  241--261 (2020)

\bibitem{scip6}
Gleixner, A., Bastubbe, M., Eifler, L., Gally, T., Gamrath, G., Gottwald, R.L.,
  Hendel, G., Hojny, C., Koch, T., L{\"u}bbecke, M., Maher, S.J., Miltenberger,
  M., M{\"u}ller, B., Pfetsch, M., Puchert, C., Rehfeldt, D., Schl{\"o}sser,
  F., Schubert, C., Serrano, F., Shinano, Y., Viernickel, J.M., Walter, M.,
  Wegscheider, F., Witt, J.T., Witzig, J.: {The SCIP Optimization Suite 6.0}.
  Tech. Rep. 18-26, ZIB, Takustr. 7, 14195 Berlin (2018)

\bibitem{miplib}
Gleixner, A., Hendel, G., Gamrath, G., Achterberg, T., Bastubbe, M., Berthold,
  T., Christophel, P., Jarck, K., Koch, T., Linderoth, J., et~al.: {MIPLIB}
  2017: data-driven compilation of the 6th mixed-integer programming library.
  Mathematical Programming Computation pp. 1--48 (2021)

\bibitem{letter-value-plot}
Hofmann, H., Kafadar, K., Wickham, H.: Letter-value plots: Boxplots for large
  data. Tech. rep., had.co.nz (2011)

\bibitem{huang2021learning}
Huang, Z., Wang, K., Liu, F., Zhen, H.l., Zhang, W., Yuan, M., Hao, J., Yu, Y.,
  Wang, J.: Learning to select cuts for efficient mixed-integer programming.
  arXiv preprint arXiv:2105.13645  (2021)

\bibitem{latticepoints}
Lodi, A., Pesant, G., Rousseau, L.M.: On counting lattice points and
  chv{\'a}tal-gomory cutting planes. In: International Conference on AI and OR
  Techniques in Constriant Programming for Combinatorial Optimization Problems.
  pp. 131--136. Springer (2011)

\bibitem{pyscipopt}
Maher, S., Miltenberger, M., Pedroso, J.P., Rehfeldt, D., Schwarz, R., Serrano,
  F.: {PySCIPOpt: Mathematical programming in python with the SCIP optimization
  suite}. In: International Congress on Mathematical Software. pp. 301--307.
  Springer (2016)

\bibitem{naoum2014recursive}
Naoum-Sawaya, J.: Recursive central rounding for mixed integer programs.
  Computers \& operations research  \textbf{43},  191--200 (2014)

\bibitem{paulus2022learning}
Paulus, M.B., Zarpellon, G., Krause, A., Charlin, L., Maddison, C.: Learning to
  cut by looking ahead: Cutting plane selection via imitation learning. In:
  International Conference on Machine Learning. pp. 17584--17600. PMLR (2022)

\bibitem{scikit-learn}
Pedregosa, F., Varoquaux, G., Gramfort, A., Michel, V., Thirion, B., Grisel,
  O., Blondel, M., Prettenhofer, P., Weiss, R., Dubourg, V., Vanderplas, J.,
  Passos, A., Cournapeau, D., Brucher, M., Perrot, M., Duchesnay, E.:
  Scikit-learn: Machine learning in {P}ython. Journal of Machine Learning
  Research  \textbf{12},  2825--2830 (2011)

\bibitem{depth}
Poirrier, L., Yu, J.: On the depth of cutting planes. arXiv preprint
  arXiv:1903.05304  (2019)

\bibitem{smola2004tutorial}
Smola, A.J., Sch{\"o}lkopf, B.: A tutorial on support vector regression.
  Statistics and computing  \textbf{14}(3),  199--222 (2004)

\bibitem{tang2020reinforcement}
Tang, Y., Agrawal, S., Faenza, Y.: {Reinforcement learning for integer
  programming: Learning to cut}. In: International Conference on Machine
  Learning. pp. 9367--9376. PMLR (2020)

\bibitem{turner_mark_2022_7433672}
Turner, M., Berthold, T., Besançon, M., Koch, T.: {Cutting Plane Selection
  with Analytic Centers and Multiregression - Density Filtering Plots} (Dec
  2022). \doi{10.5281/zenodo.7433671},
  \url{https://doi.org/10.5281/zenodo.7433671}

\bibitem{turner2022adaptive}
Turner, M., Koch, T., Serrano, F., Winkler, M.: {Adaptive Cut Selection in
  Mixed-Integer Linear Programming}. arXiv preprint arXiv:2202.10962  (2022)

\bibitem{sparsity}
Walter, M.: Sparsity of lift-and-project cutting planes. In: Operations
  Research Proceedings 2012, pp. 9--14. Springer (2014)

\bibitem{wesselmann2012implementing}
Wesselmann, F., Stuhl, U.: Implementing cutting plane management and selection
  techniques. Tech. rep., Technical report, University of Paderborn (2012)

\bibitem{wolsey2020integer}
Wolsey, L.A.: Integer programming. John Wiley \& Sons (2020)

\bibitem{lexicographicgomory}
Zanette, A., Fischetti, M., Balas, E.: Can pure cutting plane algorithms work?
  In: International Conference on Integer Programming and Combinatorial
  Optimization. pp. 416--434. Springer (2008)

\end{thebibliography}

\ifarxiv
\newpage
\appendix
\section{Cut density filtering}\label{sec:density}

In this section, we focus on another key criterion to evaluate and select cuts: cut \emph{density} or its complement \emph{sparsity}, i.e.~the number of non-zero coefficients of the cut.
In complement to the evaluation of distance measures, we therefore assess the effect of filtering cuts based on their density,
with the goal of shedding some light on a second aspect of cuts beyond the extent to which they cut off parts of the relaxation.

\subsection{Motivation}

The density of cuts is known to increase with the rank at which they are generated~\cite{sparsity}, where the rank of a cut is the number of separation rounds that were used to obtain the cut; such higher-rank cuts are known to be numerically unstable~\cite{dey2018theoretical}.
Dense cuts also typically slow down individual LP solves \cite{sparsity,chmiela2022implementation}, largely due to LU factorisations being more difficult on dense matrices.

Density is a more nuanced signal than distance measures for whether to favour a cut. Unlike the latter, one cannot consider that higher sparsity levels are systematically better.
Indeed, even though sparser cuts are numerically preferable, they can, all other things equal, be weaker by design, since one can always construct weaker sparse cuts out of dense ones.

We denote the density by \noargdensity, and relative density by \noargrelativedensity, and define them as:
\begin{align*}
    \density{\coefficients} &:= | \{i \in \{1, ..., n \} \; | \; \alpha_{i} \neq 0\} | \\
    \relativedensity{\coefficients} &:= \frac{\density{\coefficients}}{n}
\end{align*}

To the best of our knowledge, the literature is relatively sparse on the effects of density on the performance and stability of MIP solvers. Specifically for stability, consider the following cut of dimension $N \in \naturals$, where each variable $x_{i}$ is binary:
\begin{align*}
    \sum_{i \in N} x_{i} \geq 1
\end{align*}
When $N$ is large enough, the individual variables can take value $\frac{1}{N} \approx 0$ in the continuous relaxation, thereby satisfying integrality tolerances, and the cut itself.

\subsection{Computational experiments}\label{subsec:density_filtering_root}

In order to assess the impact of dense cuts on the solving process, we introduce filters that remove any cut over a given relative density threshold at the start of the cut selection algorithm, and before scoring.
We restrict ourselves to scoring cuts exclusively with efficacy for this experiment, and use relative density thresholds 5\%, 10\%, 20\%, 40\%, and 80\%, with the scoring functions denoted as \texttt{eff-05}, \texttt{eff-10}, \texttt{eff-20}, \texttt{eff-40}, and \texttt{eff-80}.

We display the results of our density experiments at the root node in Table~\ref{tab:density_root}. We see a common and expected trend of larger primal-dual differences at the end of the root node when a more restrictive density filtering is used; removing the denser cuts weakens the root-node relaxation. A corresponding decrease in solve time also follows when the filtering becomes more restrictive, although we note that this decrease is more associated with the like decrease in cuts added and amount of separation rounds, as opposed to longer LP solves. For \texttt{eff-05}, we observe that the bound performance is better, albeit still worse than \noargeff, in the $[0.8,1]$ instance subset compared to that of all other instance subsets. This is due to a few instance outliers that have substantially worse performance under \texttt{eff-05} and maximum relative densities of $0.5$ under \noargeff.
The detailed plots for this analysis are available on the archive \cite{turner_mark_2022_7433672}.

Time-specific results of our density experiments when expanded to branch-and-bound are displayed in Table~\ref{tab:time_tree}. We observe that solve time decreases for instances that would not add cuts of relative density greater than $0.4$ when using \texttt{eff-10} and \texttt{eff-20}, and decreases across all instance sets when using \texttt{eff-80}. For the most aggressive filtering \texttt{eff-05}, however, solve time increases across all instance subsets and fewer instances are solved overall.
Aside from \texttt{eff-05}, and \texttt{eff-10} for the instance subset $[0.4,1]$, all filtering methods result in at least as many instances being solved to optimality as under \noargeff. An increase in LP iteration throughput can be observed across all filtering methods and instance subsets, supporting the claim that dense cuts slow down LP solves.

Finally, we display node-specific results of our density experiments when expanded to branch-and-bound in Table~\ref{tab:nodes_tree}. We observe a slight increase of nodes needed to prove optimality when using density-based filtering, although this is partially driven by outliers as evidenced in \cite{turner_mark_2022_7433672}, highlighting the issues of performance variability for cut filtering. The exception of this node increase is \texttt{eff-10} and subsets containing instances with maximum relative densities less than $0.2$ for which density filtering improves performance. We also observe for all filtering methods and instance subsets, an increase in node throughput and LP iterations per node, empirically demonstrating that dense cuts slow down LP solves.
\begin{table}
    \begin{subtable}{1\textwidth}
    \scriptsize
    \centering
    \resizebox{\columnwidth}{!}{%
        \begin{tabular}{lcccccccccccccccccccccccccc}
            \hline
            \hline
            max-& & & \multicolumn{4}{c}{\texttt{eff-05}} && 
            \multicolumn{4}{c}{\texttt{eff-10}} && 
            \multicolumn{4}{c}{\texttt{eff-20}} && 
            \multicolumn{4}{c}{\texttt{eff-40}} && 
            \multicolumn{4}{c}{\texttt{eff-80}}\\
            \cline{4-7} \cline{9-12} \cline{14-17} \cline{19-22} \cline{24-27}
            density & instances && time & bound & cuts & round && time & bound & cuts & round && time & bound & cuts & round && time & bound & cuts & round && time & bound & cuts & round \\
            \hline
            $[0, 1]$ & 162 && 0.69 & 1.13 & 0.53 & 0.6 && 0.75 & 1.04 & 0.62 & 0.67 && 0.81 & 1.03 & 0.69 & 0.72 && 0.88 & 1.01 & 0.77 & 0.79 && 0.97 & 1.01 & 0.89 & 0.9 \\
            $[0.05, 1]$ & 134 && 0.62 & 1.16 & 0.46 & 0.54 && 0.7 & 1.05 & 0.56 & 0.61 && 0.76 & 1.03 & 0.63 & 0.67 && 0.85 & 1.02 & 0.73 & 0.75 && 0.96 & 1.01 & 0.87 & 0.88 \\
            $[0.1, 1]$ & 112 && 0.56 & 1.2 & 0.4 & 0.48 && 0.63 & 1.07 & 0.5 & 0.56 && 0.71 & 1.04 & 0.58 & 0.62 && 0.82 & 1.02 & 0.68 & 0.71 && 0.95 & 1.01 & 0.85 & 0.86 \\
            $[0.2, 1]$ & 83 && 0.47 & 1.24 & 0.32 & 0.4 && 0.54 & 1.08 & 0.41 & 0.47 && 0.62 & 1.05 & 0.47 & 0.52 && 0.76 & 1.03 & 0.59 & 0.63 && 0.93 & 1.01 & 0.8 & 0.82 \\
            $[0.4, 1]$ & 57 && 0.39 & 1.31 & 0.25 & 0.33 && 0.46 & 1.09 & 0.33 & 0.38 && 0.53 & 1.06 & 0.37 & 0.42 && 0.65 & 1.04 & 0.46 & 0.51 && 0.9 & 1.01 & 0.72 & 0.74 \\
            $[0.8, 1]$ & 23 && 0.44 & 1.12 & 0.16 & 0.24 && 0.5 & 1.09 & 0.19 & 0.26 && 0.55 & 1.08 & 0.22 & 0.29 && 0.59 & 1.07 & 0.26 & 0.33 && 0.75 & 1.04 & 0.42 & 0.48 \\
            \hline
            \hline
        \end{tabular}
    }
    \caption{Root node results. Time (s) is that to the end of all separation round. Bound is the primal-dual difference after the separation rounds. Cuts is the number of cuts added after all separation rounds. Rounds are the number of rounds.}
    \label{tab:density_root}
    \end{subtable}
    \begin{subtable}{1\textwidth}
    \scriptsize
    \centering
    \resizebox{\columnwidth}{!}{%
        \begin{tabular}{lcccccccccccccccccccccccccc}
            \hline
            \hline
            max-& & & \multicolumn{4}{c}{\texttt{eff-05}} && 
            \multicolumn{4}{c}{\texttt{eff-10}} && 
            \multicolumn{4}{c}{\texttt{eff-20}} && 
            \multicolumn{4}{c}{\texttt{eff-40}} && 
            \multicolumn{4}{c}{\texttt{eff-80}}\\
            \cline{4-7} \cline{9-12} \cline{14-17} \cline{19-22} \cline{24-27}
            density & instances && time & it & it/s & $\Delta$-sol && time & it & it/s & $\Delta$-sol && time & it & it/s & $\Delta$-sol && time & it & it/s & $\Delta$-sol && time & it & it/s & $\Delta$-sol \\
            \hline
            $[0, 1]$ & 162 && 1.15 & 1.16 & 1.04 & -20 && 0.97 & 0.99 & 1.05 & +0 && 0.99 & 1.01 & 1.03 & +2 && 1.01 & 1.01 & 1.02 & +2 && 0.99 & 0.99 & 1.01 & +3 \\
            $[0.05, 1]$ & 134 && 1.19 & 1.19 & 1.04 & -19 && 0.96 & 0.99 & 1.06 & +0 && 0.99 & 1.02 & 1.04 & +2 && 1.01 & 1.02 & 1.02 & +2 && 0.99 & 0.99 & 1.01 & +3 \\
            $[0.1, 1]$ & 112 && 1.24 & 1.26 & 1.05 & -18 && 0.95 & 0.99 & 1.08 & +0 && 0.99 & 1.02 & 1.05 & +2 && 1.01 & 1.02 & 1.03 & +2 && 0.99 & 0.99 & 1.01 & +3 \\
            $[0.2, 1]$ & 83 && 1.39 & 1.42 & 1.06 & -18 && 0.98 & 1.04 & 1.11 & +0 && 0.98 & 1.02 & 1.07 & +2 && 1.01 & 1.02 & 1.04 & +2 && 0.98 & 0.98 & 1.02 & +3 \\
            $[0.4, 1]$ & 57 && 1.63 & 1.62 & 1.05 & -14 && 1.03 & 1.09 & 1.12 & -1 && 0.99 & 1.03 & 1.07 & +0 && 1.02 & 1.04 & 1.06 & +2 && 0.98 & 0.97 & 1.02 & +3 \\
            $[0.8, 1]$ & 23 && 1.11 & 1.14 & 1.11 & +0 && 1.15 & 1.19 & 1.11 & +0 && 1.17 & 1.17 & 1.06 & +0 && 1.1 & 1.13 & 1.1 & +3 && 0.95 & 0.93 & 1.05 & +3 \\
            \hline
            \hline
        \end{tabular}
    }
    \caption{Tree results for runtime (s) to optimality or 7200s time limit. It is the number of LP iterations, and it/s the iterations per second. $\Delta$-sol is the difference in number of instance-seed pairs solved compared to \noargeff.}
    \label{tab:time_tree}
    \end{subtable}
    \begin{subtable}{1\textwidth}
    \scriptsize
    \centering
    \resizebox{\columnwidth}{!}{%
        \begin{tabular}{lcccccccccccccccccccccc}
            \hline
            \hline
            max-& & & \multicolumn{3}{c}{\texttt{eff-05}} && 
            \multicolumn{3}{c}{\texttt{eff-10}} && 
            \multicolumn{3}{c}{\texttt{eff-20}} && 
            \multicolumn{3}{c}{\texttt{eff-40}} && 
            \multicolumn{3}{c}{\texttt{eff-80}}\\
            \cline{4-6} \cline{8-10} \cline{12-14} \cline{16-18} \cline{20-22}
            density & instances && nodes & nodes/s & it/nodes && nodes & nodes/s & it/nodes && nodes & nodes/s & it/nodes && nodes & nodes/s & it/nodes && nodes & nodes/s & it/nodes \\
            \hline
            $[0, 1]$ & 134 && 1.04 & 1.17 & 0.92 && 0.97 & 1.1 & 0.97 && 1.02 & 1.07 & 0.98 && 1.03 & 1.04 & 1.0 && 1.0 & 1.01 & 1.0 \\
            $[0.05, 1]$ & 110 && 1.06 & 1.21 & 0.9 && 0.96 & 1.12 & 0.96 && 1.03 & 1.09 & 0.98 && 1.04 & 1.05 & 0.99 && 1.0 & 1.02 & 1.0 \\
            $[0.1, 1]$ & 93 && 1.11 & 1.28 & 0.86 && 0.95 & 1.15 & 0.96 && 1.03 & 1.1 & 0.97 && 1.05 & 1.06 & 0.99 && 1.0 & 1.02 & 0.99 \\
            $[0.2, 1]$ & 67 && 1.09 & 1.31 & 0.87 && 1.0 & 1.2 & 0.95 && 1.04 & 1.15 & 0.95 && 1.07 & 1.08 & 0.99 && 1.0 & 1.03 & 0.99 \\
            $[0.4, 1]$ & 47 && 1.08 & 1.32 & 0.88 && 1.01 & 1.22 & 0.96 && 1.03 & 1.16 & 0.97 && 1.1 & 1.12 & 0.99 && 0.99 & 1.04 & 0.99 \\
            $[0.8, 1]$ & 21 && 1.1 & 1.38 & 0.88 && 1.1 & 1.31 & 0.91 && 1.08 & 1.25 & 0.95 && 1.08 & 1.22 & 0.94 && 0.99 & 1.1 & 0.97 \\
            \hline
            \hline
        \end{tabular}
    }
    \caption{Tree results. Nodes is the number of nodes to prove optimality. Nodes/s is the node throughput per second, and it/nodes is the number of LP iterations per node.}
    \label{tab:nodes_tree}
    \end{subtable}
    \caption{Summary of density results. Rows labelled $[d,1]$ consider instances where the default runs with efficacy added a cut of relative density of at least $d$. Entries are the instance-wise shifted geometric mean normalised by the shifted geometric mean using \noargeff. The shifts are: (time, 1), (bound, 1), (cuts, 10), (round, 1), (it, 100), (it/s, 10), (nodes, 10), (nodes/s, 1), (it/n, 10).}
    \label{tab:density}
\end{table}


\else
\fi
\end{document}